\documentclass[11pt]{article}
\usepackage{graphicx}
\usepackage{latexsym}
\usepackage{amssymb}

\def\beq{\begin{eqnarray}}
\def\eeq{\end{eqnarray}}

\def\geo{^{\mbox{\tiny geo}}}
\def\opt{^{\mbox{\tiny opt}}}
\def\rhat{\hat r}
\def\rhatrev{\underline{\hat r}}
\def\rhatSIS{\hat r_{\mbox{\tiny SIS}}}
\def\rhatbridge{\hat r_{\mbox{\tiny bridge}}}
\def\rhatAIS{\hat r_{\mbox{\tiny AIS}}}
\def\rhatLIS{\hat r_{\mbox{\tiny LIS}}}
\def\rhatLISave{\hat r_{\mbox{\tiny LIS-ave}}}
\def\rhatLISrev{\underline{\hat r}_{\,\mbox{\tiny LIS}}}
\def\rhatLISbridged{\hat r_{\mbox{\tiny LIS-bridged}}}
\def\Var{\mbox{Var}}
\def\Cor{\mbox{Cor}}

\begin{document}

\fontsize{11}{16pt}\selectfont

\begin{center}

{\small Technical Report No.\ 0511,
 Department of Statistics, University of Toronto}

\vspace*{0.45in}

{\LARGE \bf Estimating Ratios of Normalizing Constants Using \\[6pt]
            Linked Importance Sampling}

\vspace*{9pt}

{\large Radford M. Neal}\\[4pt]
 Department of Statistics and Department of Computer Science \\
 University of Toronto, Toronto, Ontario, Canada \\
 \texttt{http://www.cs.utoronto.ca/$\sim$radford/} \\
 \texttt{radford@stat.utoronto.ca}\\[6pt]

 8 November 2005

\end{center}

\vspace*{8pt}

\noindent \textbf{Abstract.}\ \ Ratios of normalizing constants for
two distributions are needed in both Bayesian statistics, where they
are used to compare models, and in statistical physics, where they
correspond to differences in free energy.  Two approaches have long
been used to estimate ratios of normalizing constants.  The `simple
importance sampling' (SIS) or `free energy perturbation' method uses a
sample drawn from just one of the two distributions.  The `bridge
sampling' or `acceptance ratio' estimate can be viewed as the ratio of
two SIS estimates involving a bridge distribution.  For both methods,
difficult problems must be handled by introducing a sequence of
intermediate distributions linking the two distributions of interest,
with the final ratio of normalizing constants being estimated by the
product of estimates of ratios for adjacent distributions in this
sequence.  Recently, work by Jarzynski, and independently by Neal, has
shown how one can view such a product of estimates, each based on
simple importance sampling using a single point, as an SIS estimate on
an extended state space.  This `Annealed Importance Sampling' (AIS)
method produces an exactly unbiased estimate for the ratio of
normalizing constants even when the Markov transitions used do not
reach equilibrium.  In this paper, I show how a corresponding `Linked
Importance Sampling' (LIS) method can be constructed in which the
estimates for individual ratios are similar to bridge sampling
estimates.  As a further elaboration, bridge sampling rather than
simple importance sampling can be employed at the top level for both
AIS and LIS, which sometimes produces further improvement.  I show
empirically that for some problems, LIS estimates are much more
accurate than AIS estimates found using the same computation time,
although for other problems the two methods have similar performance.
Like AIS, LIS can also produce estimates for expectations, even when
the distribution contains multiple isolated modes.  AIS is related to
the `tempered transition' method for handling isolated modes, and to a
method for `dragging' fast variables.  Linked sampling methods similar
to LIS can be constructed that are analogous to tempered transitions
and to this method for dragging fast variables, which may sometimes
work better than those analogous to AIS.

\newpage

\section{\hspace*{-7pt}Introduction}\label{sec-intro}\vspace*{-10pt}

Consider two distributions on the same space, with probability mass or
density functions $\pi_0(x) = p_0(x)/Z_0$ and $\pi_1(x) =
p_1(x)/Z_1$.  Suppose that we are not able to directly compute $\pi_0$
and $\pi_1$, but only $p_0$ and $p_1$, since we do not know the
normalizing constants, $Z_0$ and $Z_1$.  We wish to find a Monte Carlo
estimate for the ratio of these normalizing constants, $Z_1/Z_0$,
which we sometimes denote by $r$, using samples of values drawn (at
least approximately) from $\pi_0$ and from $\pi_1$.  Sometimes, we may
know $Z_0$, in which case we can arrange for it to be one, so that
estimation of this ratio will give the numerical value of $Z_1$.
Other times, we will be able to obtain only the ratio of normalizing
constants, but this may be sufficient for our purposes.

In statistical physics, $x$ represents the state of some physical
system, and the distributions are typically `canonical' distributions
having the following form (for $j=0,1$):
\beq
  p_j(x) & = & \exp(-\beta_j U(x,\lambda_j))
\label{eq-canonical}
\eeq
where $U(x,\lambda_j)$ is an `energy' function, which may depend on the
parameter $\lambda_j$, and $\beta_j$ is the inverse temperature of
system $j$.  Many interesting properties of the systems are related
to the `free energy', defined as $-\log(Z_j)\,/\,\beta_j$.  Often, only
the difference in free energy between systems $0$ and $1$ is relevant, 
and this is determined by the ratio $Z_1/Z_0$.

In Bayesian statistics, $x$ comprises the parameters and latent
variables for some statistical model, $\pi_0$ is the prior
distribution for these quantities (for which the normalizing constant
is usually known), and $\pi_1$ is the posterior distribution given the
observed data.  We can compute $p_1(x)$ as the product of the prior
density for $x$ and the probability of the data given $x$, but the
normalizing constant, $Z_1$, is difficult to compute.  We can
interpret $Z_1$ as the `marginal likelihood' --- the probability of
the observed data under this model, integrating over possible values
of the model's parameters and latent variables.  The marginal
likelihood for a model indicates how well it is supported
by the data.

Although I will use simple distributions as illustrations in this
paper, in real applications, $x$ is usually high dimensional, and at
least one of $\pi_0$ and $\pi_1$ is usually quite complex.
Accordingly, sampling from these distributions generally requires use
of Markov chain methods, such as the venerable Metropolis algorithm
(Metropolis, \textit{et al} 1953).  See (Neal 1993) for a review of
Markov chain sampling methods.  Sometimes, however, $\pi_0$ will be
relatively simple, and independent points drawn from it can be
generated efficiently, as would often be the case with the prior
distribution for a Bayesian model, or for a physical system at
infinite temperature ($\beta_0=0$).

Many methods for estimating ratios of normalizing constants from Monte
Carlo data have been investigated in the physics literature (for a
review, see (Neal 1993, Section 6.2)), and later rediscovered in the
statistics literature (Gelman and Meng 1998).  A logical method to
start with is `simple importance sampling' (SIS), also called `free energy
perturbation', based on the following identity, which can 
easily be proved on the assumption that no region having zero probability 
under $\pi_0$ has non-zero probability under $\pi_1$:
\beq
  {Z_1 \over Z_0} & = & E_{\pi_0}\! \left[ {p_1(X) \over p_0(X)} \right]
   \ \ \approx \ \ {1 \over N} \sum_{i=1}^N {p_1(x^{(i)}) \over p_0(x^{(i)})}
   \ \ =\ \  {1 \over N} \sum_{i=1}^N \rhatSIS^{(i)}
   \ \ =\ \ \rhatSIS
\label{eq-simple} 
\eeq 
In the above equation, $E_{\pi_0}$ denotes an expectation with 
respect to the distribution
$\pi_0$, which is estimated by a Monte Carlo average over points
$x^{(i)},\ldots,x^{(N)}$ drawn from $\pi_0$ (either independently, or using a
Markov chain sampler).  
Here and later, $\hat r_{\mbox{\tiny M}}$ will denote an estimate of 
$r=Z_1/Z_0$, found by method M.  If this estimate is an average of 
unbiased estimates based on a number of samples, these individual 
estimates will be denoted by $\hat r_{\mbox{\tiny M}}^{(i)}$.  

The simple importance sampling estimate, $\rhatSIS$, will be poor
if $\pi_0$ and $\pi_1$ are not close enough --- in particular, if any
region with non-negligible probability under $\pi_1$ has very small
probability under $\pi_0$.  Such a region would have an important
effect on the value of $r$, but very little information about it would
be contained in the sample from $\pi_0$.  In such a situation, it may
be possible to obtain a good estimate by introducing intermediate
distributions.  Parameterizing these distributions in some way using
$\eta$, we can define a sequence of distributions,
$\pi_{\eta_0},\ldots,\pi_{\eta_n}$, with $\eta_0=0$ and $\eta_n=1$ so
that the first and last distributions in the sequence are $\pi_0$ and
$\pi_1$, with the intermediate distributions interpolating between
them.  We can then write
\beq
  {Z_1 \over Z_0} & = & \prod_{j=0}^{n-1} {Z_{\eta_{j+1}} \over Z_{\eta_j}}
\label{eq-intermed}
\eeq
Provided that $\pi_{\eta_{j+1}}$ and $\pi_{\eta_j}$ are close enough, 
we can estimate each of the factors 
$Z_{\eta_{j+1}}/Z_{\eta_j}$ using simple
importance sampling, and from these estimates obtain an estimate for $Z_1/Z_0$.

We can obtain good estimates in a wider range of situations, or using
fewer intermediate distributions (sometimes none), by applying a
technique introduced by Bennett (1976), who called it the `acceptance
ratio' method.  This method was later rediscovered by Meng and Wong
(1996), who called it `bridge sampling'.  Lu, Singh, and Kofke (2003)
provide a recent review and assessment.  One way of viewing this
method is that it replaces the simple importance sampling estimate for
$Z_1/Z_0$ by a ratio of estimates for $Z_*/Z_0$ and $Z_*/Z_1$, where
$Z_*$ is the normalizing constant for a `bridge distribution',
$\pi_*(x) = p_*(x)/Z_*$, which is chosen so that it is overlapped by
both $\pi_0$ and $\pi_1$.  Using simple importance sampling estimates
for $Z_*/Z_0$ and $Z_*/Z_1$, we can obtain the estimate
\beq
  {Z_1 \over Z_0} & = & 
     E_{\pi_0}\! \left[ {p_*(X) \over p_0(X)} \right] \, \Big/\,
     E_{\pi_1}\! \left[ {p_*(X) \over p_1(X)} \right]
  \ \ \approx \ \ 
 {1 \over N_0} \sum_{k=1}^{N_0} {p_*(x_{0,k}) \over p_0(x_{0,k})} \ \Big/\
 {1 \over N_1} \sum_{k=1}^{N_1} {p_*(x_{1,k}) \over p_1(x_{1,k})} 
   \ \ =\ \ \rhatbridge\ \ \ \ \ 
\label{eq-bridge}
\eeq
where $x_{0,1},\ldots,x_{0,N_0}$ are drawn from $\pi_0$ and
$x_{1,1},\ldots,x_{1,N_1}$ are drawn from $\pi_1$.  

One simple choice for the bridge distribution is the `geometric' bridge:
\beq
  p\geo_*(x) & = & \sqrt{p_0(x)p_1(x)}
\label{eq-geo-bridge}
\eeq
which is in a sense half-way between $\pi_0$ and $\pi_1$.  
As discussed by Bennett (1976) and by Meng and Wong (1996), the asymptotically
optimal choice of bridge distribution is
\beq
  p\opt_*(x) & = & { p_0(x)p_1(x) \over r (N_0/N_1) p_0(x)\, +\, p_1(x)}
\label{eq-opt-bridge}
\eeq
where $r=Z_1/Z_0$.  Of course, we cannot use this bridge distribution
in practice, since we do not know $r$.  We can use a preliminary guess
at $r$ to define an initial bridge distribution, however, which will
give us a bridge sampling estimate for $Z_1/Z_0$.  Using this estimate
as the new value of $r$, we can refine our bridge distribution, iterating
this process as many times as desired.  The result of this iteration can
also be viewed as a maximum likelihood estimate for $r$, as discussed by
Shirts, \textit{et~al} (2003), who argues on this basis that it is
asymptotically as good as any estimate for $r$.   I have found that
estimates with $r$ set iteratively are often better than those found
with the true value of $r$ (which does not contradict optimality of the true
value for a fixed choice of bridge distribution).

If $\pi_0$ and $\pi_1$ do not overlap sufficiently, no bridge
distribution will produce good estimates, and we will have to
introduce intermediate distributions as in
equation~(\ref{eq-intermed}).  Note, however, that the bridge sampling
estimate with either of the above bridge distributions converges
to the correct ratio asymptotically as long there is some region that
has non-zero probability under both $\pi_0$ and $\pi_1$, a much weaker
requirement than that for simple importance sampling.

This advantage of bridge sampling over SIS can be seen in a simple
example involving distributions that are uniform over an interval of the
reals.  Let $p_0(x) = I_{(0,3)}(x)$ and $p_1(x)=I_{(2,4)}(x)$, so that
$Z_0=3$ and $Z_1=2$.  The simple importance sampling estimate of
equation~(\ref{eq-simple}) does not work, as it converges to $1/3$
rather than $2/3$.  However, using a bridge distribution with
$p_*(x)=I_{(2,3)}$, which is effectively what both $p_*\opt$ and
$p_*\geo$ will be in this example, the bridge sampling estimate of
equation~(\ref{eq-bridge}) converges to the correct value, since the
numerator converges to $1/3$ and the denominator to $1/2$.

Although both simple importance sampling and bridge sampling have been
successfully used in many applications, they have some deficiencies.
One issue is that although the SIS estimate of
equation~(\ref{eq-simple}) is unbiased for $Z_1/Z_0$, the bridge
sampling estimate of equation~(\ref{eq-bridge}) is not, and the same
would appear to be the case for an estimate using intermediate
distributions (via equation~(\ref{eq-intermed})).  This is of no
direct importance, particularly since we are often more interested in
$\log(Z_1/Z_0)$ than in $Z_1/Z_0$ itself.  However, it does preclude
averaging independent replications of the bridge sampling estimate to
obtain a better estimate, since the bias would prevent convergence to
the correct value as the number of replications increases.  A more
vexing difficulty is that, except sometimes for $\pi_0$, sampling from
the distributions $\pi_{\eta}$ must usually be done by Markov chain
methods, which approach the desired distribution only asymptotically.
To speed convergence, the Markov chain for sampling $\pi_{\eta_j}$ is
often started from the last state sampled for $\pi_{\eta_{j-1}}$, but
it is unclear how many iterations should then be discarded before an
adequate approximation to the correct distribution is reached.

Surprisingly, these difficulties can be completely overcome when using
simple importance sampling with a single point.  As shown by Jarzynski
(1997, 2001), and later independently by myself (Neal 2001), an estimate for
$Z_1/Z_0$ using intermediate distributions as in
equation~(\ref{eq-intermed}) will be exactly unbiased if each of the
ratios $Z_{\eta_{j+1}}/Z_{\eta_j}$ is estimated using the simple
importance sampling estimate of equation~(\ref{eq-simple}) with $N=1$,
sampling each distribution with a Markov chain update starting with the
point for the previous distribution.
Averaging the estimates obtained from $M$ independent replications of this 
process (called `runs') produces the following estimate:
\beq
  {Z_1 \over Z_0} & \approx & 
    {1 \over M}\, \sum_{i=1}^M\, \prod_{j=0}^{n-1}\,
    {p_{\eta_{j+1}}(x^{(i)}_j) \over p_{\eta_j}(x^{(i)}_j)}
    \ \ =\ \ {1 \over M} \sum_{i=1}^M \rhatAIS^{(i)}
    \ \ =\ \ \rhatAIS
\label{eq-ais-est}
\eeq
Here, $x^{(1)}_0,\ldots,x^{(M)}_0$ are drawn independently from $\pi_0$, 
and each $x^{(i)}_j$ for $j>0$ is generated by applying a Markov chain 
transition that leaves $\pi_{\eta_j}$ invariant to $x^{(i)}_{j-1}$.  This 
single Markov transition (which could, however, consist of several Metropolis 
or other updates if we so choose), will usually not be enough to reach 
equilibrium, but the estimate $\rhatAIS$ is nevertheless exactly unbiased, and 
will converge to the true value as $M$ increases, provided that no region
having zero probability under $\pi_{\eta_j}$ has non-zero probability
under $\pi_{\eta_{j+1}}$.  This can be proved by showing how the
estimate above can be seen as a simple importance sampling estimate on an 
extended state space that includes the values sampled for the intermediate 
distributions.

I call this method `Annealed Importance Sampling' (AIS), since the
sequence of distributions used often corresponds to an `annealing'
procedure, in which the temperature is gradually decreased.  As I
discuss in (Neal 2001), this allows the procedure to sample different
isolated modes of the distribution on different runs, properly
weighting the points obtained from each of these runs to produce the
correct probability for each mode.  AIS is related to an earlier
method for moving between isolated modes that I call `tempered
transitions' (Neal 1996).  In a recent paper (Neal 2004), I show how
tempered transitions can be modified to produce a method for efficient
Markov chain sampling when some of the state variables are `fast' ---
ie, when it is possible to more quickly recompute the probability of a
state when only these fast variables change than when the other `slow'
variables change as well.  In this method, the fast variables are
`dragged' through intermediate distributions in order to produce more
appropriate values to go with a proposed change to the slow variables.
Deciding whether to accept the final proposal involves what is in
effect an estimate of the ratio of normalizing constants for the
conditional distributions of the fast variables.

In this paper, I show how the ideas behind Annealed Importance
Sampling and bridge sampling can be combined.  I call the resulting
method `Linked Importance Sampling' (LIS), since the two samples
needed for bridge sampling are linked by a single state that is used
in both.  Intermediate distributions can be used, with each
distribution being linked by a single state to the next distribution.
In contrast to bridge sampling, LIS estimates are unbiased, and as is
the case for AIS, they remain exactly unbiased even when intermediate
distributions are used, and when sampling is done using Markov chain
transitions that have not converged to their equilibrium
distributions.

Crooks (2000) mentions a different way of combining AIS with bridge
sampling --- since AIS estimates are simple importance sampling
estimates on an extended state space, we can combine `forward' and
`reverse' estimates to produce a bridge sampling estimate that may be
superior.  I will call this method `bridged AIS'.  Similarly,
such a top-level application of bridge sampling can be combined with
the low-level application of bridge sampling in LIS, giving what I
call `bridged LIS'.

Using tests on sequences of one-dimensional distributions, I
demonstrate that for some problems LIS is much more efficient than AIS
--- a result that should be expected, since in extreme cases, such as
for the uniform distributions discussed above, the simple importance
sampling estimates underlying AIS do not converge to the correct
answer even asymptotically, whereas bridge sampling estimates do.  For
some other problems, however, AIS and LIS perform about equally well.
The bridged version of AIS sometimes performs much better than the
unbridged version, but still performs less well than LIS and its
bridged version on some problems.  I also analyse the asymptotic
properties of AIS and LIS for some types of distribution, providing
additional insight into their behaviour.

Variants of tempered transitions and of my method for dragging fast
variables can be constructed that are analogous to LIS rather than to
AIS.  I discuss the `linked' variant of tempered transitions briefly,
and include a more detailed description of a linked version of
dragging, which may sometimes be better than the version related to
AIS.  I conclude by discussing some possibilities for future research.

\section{\hspace*{-7pt}The Linked Importance Sampling 
         procedure}\label{sec-lis}\vspace*{-10pt}

Assume that we can evaluate the unnormalized probability or density
functions $p_{\eta}(x)$, for any value of the parameter $\eta$, with
the normalized form of such a distribution being denoted by
$\pi_{\eta}$.  The values $\eta=0$ and $\eta=1$ define the two
distributions we are interested in, for which the normalizing
constants are $Z_0$ and $Z_1$.  A sequence of $n\!-\!1$ intermediate
values for $\eta$ define distributions that will assist in estimating
the ratio of these normalizing constants, $r=Z_1/Z_0$.  We denote the
values of $\eta$ for the distributions used by $\eta_0,\ldots,\eta_n$,
with $\eta_0=0$ and $\eta_n=1$.  Typically, $\eta_j<\eta_{j+1}$ for
all $j$.

For problems in statistical physics, $\eta$ might be proportional to
the inverse temperature, $\beta$, of equation~(\ref{eq-canonical}), or
might map to a value for $\lambda$.  For a Bayesian inference
problem, $\eta$ might be a power that the likelihood is raised to, so
that $\eta=0$ causes the data to be ignored, and $\eta=1$ gives full
weight to the data; the ratio $Z_1/Z_0$ will then be the marginal
likelihood.  In both of these examples, progressing in small steps
from $\eta=0$ to $\eta=1$ is not only useful in estimating $Z_1/Z_0$,
but also often has an `annealing' effect, which helps avoid being
trapped in a local mode of the distribution.

\subsection{\hspace*{-4pt}Details of the LIS procedure}\vspace*{-4pt}

For each distribution, $\pi_{\eta}$, assume we have a pair of Markov chain
transition probability (or density) functions, denoted by $T_{\eta}(x,x')$ 
and $\underline{T}_{\eta}(x,x')$, satisfying $\int T_{\eta}(x,x') dx' = 1$
and $\int \underline{T}_{\eta}(x,x') dx' = 1$, for which the following mutual
reversibility relationship holds: 
\beq
  \pi_{\eta}(x)\,T_{\eta}(x,x') & = & 
     \pi_{\eta}(x')\,\underline{T}_{\eta}(x',x),\ \ \ \
  \mbox{for all $x$ and $x'$}
\label{eq-rev}
\eeq
From this relationship, one can easily show that both $T_{\eta}$ and
$\underline{T}_{\eta}$ leave $\pi_{\eta}$ invariant --- ie, that 
$\int \pi_{\eta}(x)
T_{\eta}(x,x') dx = \pi_{\eta}(x')$, and the same for $\underline{T}_{\eta}$.
If $T_{\eta}$ is reversible (ie, satisfies `detailed balance'), then 
$\underline{T}_{\eta}$ will be the same as $T_{\eta}$.  Non-reversible 
transitions often arise when components of state are updated in some 
predetermined order, in which case the reverse transition simply updates
components in the opposite order.  As a special case, $T_{\eta}$ might
draw the next state from $\pi_{\eta}$ independently of the current state.
Such independent sampling may often be possible for $T_0$.

These Markov chain transitions are used to obtain samples that are
approximately drawn from each of the $n\!+\!1$ distributions,
$\pi_{\eta_0},\ldots,\pi_{\eta_n}$.  We assume that we can begin
sampling from $\pi_0$ by drawing a single point independently from
$\pi_0$.  For $j>0$, we begin sampling from $\pi_{\eta_j}$ by
selecting a link state, $x_{j-1*j}$, from the sample associated with
$\pi_{\eta_{j-1}}$.  For all $j$, we produce a sample of $K_j\!+\!1$
states from this starting point by applying a total of $K_j$ forward
($T_{\eta_j}$) or reversed ($\underline T_{\eta_j}$) Markov
transitions.  Link states are selected using bridge distributions,
$p_{j*j+1}$, which are defined in terms of $p_{\eta_j}$ and
$p_{\eta_{j+1}}$, perhaps using the form of
equation~(\ref{eq-geo-bridge}) or~(\ref{eq-opt-bridge}), with $p_0$
replaced by $p_{\eta_j}$ and $p_1$ by $p_{\eta_{j+1}}$.

In detail, the Linked Importance Sampling procedure produces $M$ estimates,
$\rhatLIS^{(1)},\ldots,\rhatLIS^{(M)}$, that are averaged to produce
the final estimate, $\rhatLIS$.  Each $\rhatLIS^{(i)}$ is
obtained by performing the following:\vspace*{5pt}

\begin{center}\bf The LIS Procedure\end{center}\vspace*{-5pt}

\begin{figure}[t]

\centerline{\includegraphics[width=6.5in]{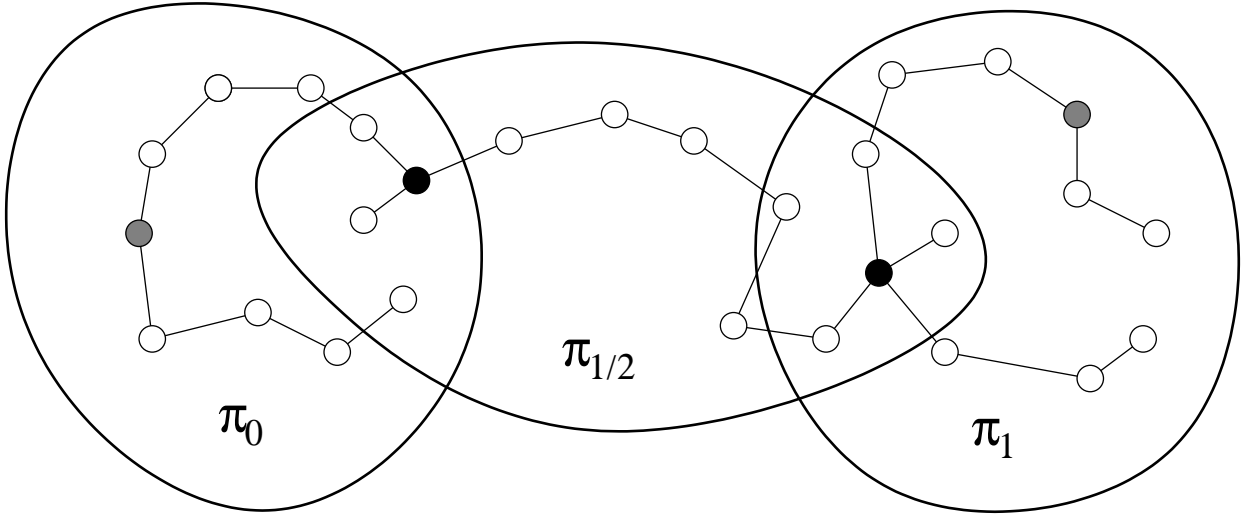}}

\caption[]{An illustration of Linked Importance Sampling.  One
intermediate distribution is used, with $\eta_1=1/2$.  The
distributions $\pi_0$, $\pi_{1/2}$, and $\pi_1$ are represented by
ovals enclosing the regions of high probability under each
distribution.  Nine Markov chain transitions are performed at each
stage.  The two link states are shown as black dots.  The initial and
final states (indexed by $\nu_0$ and $\mu_n$) are shown as gray dots.
Other states generated by the forward and reverse Markov chain
transitions are shown as empty dots.  For this run, $\nu_0\!=\!4$,
$\mu_0\!=\!9$, $\nu_1\!=\!1$, $\mu_1\!=\!8$, $\nu_2\!=\!3$, and
$\mu_2\!=\!7$.}\label{fig-lis} 

\end{figure}

\begin{enumerate}
\item[1)] Pick an integer $\nu_0$ uniformly at random from $\{0,\ldots,K_0\}$,
          and then set $x_{0,\nu_0}$ to a value drawn from $\pi_{\eta_0}$.
\item[2)] For $j\,=\,0,\ldots,n$, sample $K_j\!+\!1$ states drawn (at 
  least approximately) from $\pi_{\eta_j}$  as follows:
  \begin{enumerate}
  \item[a)] If $j>0$:\ \ Pick an integer $\nu_j$ uniformly at random from 
            $\{0,\ldots,K_j\}$, and then set $x_{j,\nu_j}$ to $x_{j-1*j}$.
  \item[b)] For $k\,=\,\nu_j+1,\ldots,K_j$, draw $x_{j,k}$ according to the
            forward Markov chain transition probabilities 
            $T_{\eta_j}(x_{j,k-1},x_{j,k})$.  (If $\nu_j=K_j$, do nothing in 
            this step.)
  \item[c)] For $k\,=\,\nu_j-1,\ldots,0$, draw $x_{j,k}$ according to the 
            reverse Markov chain transition probabilities
            $\underline{T}_{\eta_j}(x_{j,k+1},x_{j,k})$. (If $\nu_j=0$, do 
            nothing in this step.)
  \item[d)] If $j<n$:\ \ Pick a value for $\mu_j$ from
            $\{0,\ldots,K_j\}$ according to the following 
            probabilities:\vspace*{-2pt}
            \beq
              \Pi_0(\mu_j\,|\,x_j) & = & 
                 {p_{j*j+1}(x_{j,\mu_j}) \over p_{\eta_j}(x_{j,\mu_j})}
                 \ \Big/\
                 \sum_{k=0}^{K_j} {p_{j*j+1}(x_{j,k}) \over p_{\eta_j}(x_{j,k})}
            \label{eq-pmuj}
            \eeq
            and then set $x_{j*j+1}$ to $x_{j,\mu_j}$.
  \end{enumerate}
\item[3)] Set $\mu_n$ to a value chosen uniformly at random from 
          $\{0,\ldots,K_n\}$.  (This selection has no effect on
          the estimate, but is used in the proof of correctness.)
\item[4)] Compute the estimate from this run as follows:
\beq
   \rhatLIS^{(i)} & = & \prod_{j=0}^{n-1} \left[
     {1 \over K_j+1}\, \sum_{k=0}^{K_j}\,
          { p_{j*j+1}(x_{j,k}) \over p_{\eta_j}(x_{j,k}) }
     \ \Big/\
     {1 \over K_{j+1}+1}\, \sum_{k=0}^{K_{j+1}}\, 
          { p_{j*j+1}(x_{j+1,k}) \over p_{\eta_{j+1}}(x_{j+1,k}) }
     \right]
\label{eq-lis}
\eeq
(Note that most of the factors of $1/(K_j\!+\!1)$ and
$1/(K_{j+1}\!+\!1)$ cancel, giving a final result of 
$(K_n\!+\!1)\,/\,(K_0\!+\!1)$, but the redundant factors
are retained above for clarity of meaning.)\vspace*{-6pt}
\end{enumerate}
The result of performing steps (1) through (3) is illustrated in 
Figure~\ref{fig-lis}.  After $M$ runs of this procedure, the final
estimate is computed as
\beq
  \rhatLIS & = & {1 \over M} \sum_{i=1}^M \rhatLIS^{(i)}
\eeq

The crucial aspect of Linked Importance Sampling is that when moving
from distribution $\pi_{\eta_j}$ to $\pi_{\eta_{j+1}}$, a link state,
$x_{j*j+1}$, is randomly selected from among the sample of points
$x_{j,1},\ldots,x_{j,K_j+1}$ that are associated with $\pi_{\eta_j}$.
We can view the link state as part of the sample associated with
$\pi_{\eta_{j+1}}$ as well as that associated with $\pi_{\eta_j}$.
Accordingly, when using the `optimal' bridge of
equation~(\ref{eq-opt-bridge}), I will set $N_0/N_1$ to
$(K_j\!+\!1)/(K_{j+1}\!+\!1)$, though the proof of optimality for
bridge sampling does not guarantee that this is an optimal choice when
using this bridge distribution for LIS.

\subsection{\hspace*{-4pt}Proof that LIS estimates are unbiased}\vspace*{-4pt}

In order to prove that $\rhatLIS^{(i)}$ is an unbiased estimate of
$r=Z_1/Z_0$, we can regard steps (1) through (3) above as defining a 
distribution,
$\Pi_0$, over all the quantities involved in the procedure --- namely, 
$x_j$, $\mu_j$, and $\nu_j$, for $j=0,\ldots,n$, with $x_j$ representing
$x_{j,0},\ldots,x_{j,K_j}$.  We then 
consider the procedure for generating these same quantities in reverse,
which operates as follows:\vspace*{5pt}

\pagebreak

\begin{center}\bf The Reverse LIS Procedure\end{center}\vspace*{-5pt}

\begin{enumerate}
\item[1)] Pick an integer $\mu_n$ uniformly at random from $\{0,\ldots,K_n\}$,
          and then set $x_{n,\mu_n}$ to a value drawn from $\pi_{\eta_n}$.
\item[2)] For $j\,=\,n,\ldots,0$, sample $K_j\!+\!1$ states drawn (at 
  least approximately) from $\pi_{\eta_j}$  as follows:
  \begin{enumerate}
  \item[a)] If $j<n$:\ \ Pick an integer $\mu_j$ uniformly at random from 
            $\{0,\ldots,K_j\}$, and then set $x_{j,\mu_j}$ to $x_{j*j+1}$.
  \item[b)] For $k\,=\,\mu_j+1,\ldots,K_j$, draw $x_{j,k}$ according to the
            forward Markov chain transition probabilities 
            $T_{\eta_j}(x_{j,k-1},x_{j,k})$.  (If $\mu_j=K_j$, do nothing 
            in this step.)
  \item[c)] For $k\,=\,\mu_j-1,\ldots,0$, draw $x_{j,k}$ according to the 
            reverse Markov chain transition probabilities
            $\underline{T}_{\eta_j}(x_{j,k+1},x_{j,k})$. (If $\mu_j=0$, 
            do nothing in this step.)
  \item[d)] If $j>0$:\ \ Pick a value for $\nu_j$ from 
            $\{0,\ldots,K_j\}$ according to the following 
            probabilities:\vspace*{-3pt}
            \beq
              \Pi_1(\nu_j\,|\,x_j) & = & 
                 {p_{j-1*j}(x_{j,\nu_j}) \over p_{\eta_j}(x_{j,\nu_j})}
                 \ \Big/\
                 \sum_{k=0}^{K_j} {p_{j-1*j}(x_{j,k}) \over p_{\eta_j}(x_{j,k})}
            \label{eq-pnuj}
            \eeq
            and then set $x_{j-1*j}$ to $x_{j,\nu_j}$.
  \end{enumerate}
\item[3)] Set $\nu_0$ to a value chosen uniformly at random from 
          $\{0,\ldots,K_0\}$.\vspace*{-6pt}
\end{enumerate}
This reverse procedure also defines a distribution over all the 
quantities generated ($x_j$, $\mu_j$, and $\nu_j$ for $j=0,\ldots,n$),
which will be denoted by $\Pi_1$.

We now define the unnormalized probability (density) functions
$P_0(x,\mu,\nu) = Z_0 \Pi_0(x,\mu,\nu)$ and 
$P_1(x,\mu,\nu) = Z_1 \Pi_1(x,\mu,\nu)$.  The ratio of normalizing constants
for these distributions is obviously $r=Z_1/Z_0$.  We can estimate this
ratio by simple importance sampling, using the ratios
\beq
   {P_1(x,\mu,\nu) \over P_0(x,\mu,\nu)} & = &
   { Z_1\, \Pi_1(\mu_n)\, \pi_{\eta_n}(x_{n,\mu_n})\,
     \prod\limits_{j=0}^{n-1} \Pi_1(\mu_j)\,
     \prod\limits_{j=0}^n \Pi_1(x_j\,|\,\mu_j,x_{j,\mu_j})\,
     \prod\limits_{j=1}^{n} \Pi_1(\nu_j\,|\,x_j)\, \Pi_1(\nu_0)
     \over 
     Z_0\, \Pi_0(\nu_0)\, \pi_{\eta_0}(x_{0,\nu_0})\,
     \prod\limits_{j=1}^n \Pi_0(\nu_j)\, 
     \prod\limits_{j=0}^n \Pi_0(x_j\,|\,\nu_j,x_{j,\nu_j})\,
     \prod\limits_{j=0}^{n-1} \Pi_0(\mu_j\,|\,x_j)\, \Pi_0(\mu_n)
   }\ \ \
\label{eq-ratio01}
\eeq

From Steps (2b) and (2c) of the forward and reverse procedures, along 
with the mutual reversibility relationship of equation~(\ref{eq-rev}), we see 
that
\beq
  \Pi_0(x_j\,|\,\nu_j,x_{j,\nu_j}) 
  & = & 
  \prod_{k=\nu_j+1}^n\!\! T_{\eta_j}(x_{j,k-1},x_{j,k})\ \cdot\
  \prod_{k=0}^{\nu_j-1} \underline{T}_{\eta_j}(x_{j,k+1},x_{j,k}) \\[4pt]
  & = & 
  \prod_{k=\nu_j+1}^n\!\! T_{\eta_j}(x_{j,k-1},x_{j,k})\ \cdot\
  \prod_{k=0}^{\nu_j-1} T_{\eta_j}(x_{j,k},x_{j,k+1})\,
                        {\pi_{\eta_j}(x_{j,k})\over\pi_{\eta_j}(x_{j,k+1})}
  \\[4pt]
  & = & 
  {\pi_{\eta_j}(x_{j,0})\over\pi_{\eta_j}(x_{j,\nu_j})}\
  \prod_{k=1}^n\, T_{\eta_j}(x_{j,k-1},x_{j,k})
\label{eq-chain1}
\eeq
and similarly,
\beq
  \Pi_1(x_j\,|\,\mu_j,x_{j,\mu_j}) 
  & = & 
  {\pi_{\eta_j}(x_{j,0})\over\pi_{\eta_j}(x_{j,\mu_j})}\
  \prod_{k=1}^n\, T_{\eta_j}(x_{j,k-1},x_{j,k})
\label{eq-chain2}
\eeq
From this, we see that parts of the ratio in equation~(\ref{eq-ratio01})
can be written as
\beq
   { Z_1\,\pi_{\eta_n}(x_{n,\mu_n})\, 
     \prod\limits_{j=0}^n \, \Pi_1(x_j\,|\,\mu_j,x_{j,\mu_j})\,
     \over 
     Z_0\,\pi_{\eta_0}(x_{0,\nu_0})\,
     \prod\limits_{j=0}^n \, \Pi_0(x_j\,|\,\nu_j,x_{j,\nu_j})\,
   } 
   & = & 
   {p_{\eta_n}(x_{n,\mu_n}) \over p_{\eta_0}(x_{0,\nu_0})}\,
   \prod_{j=0}^n\, {\pi_{\eta_j}(x_{j,\nu_j}) \over \pi_{\eta_j}(x_{j,\mu_j})}
   \ \ =\ \
   \prod_{j=0}^{n-1}\, 
      {p_{\eta_{j+1}}(x_{j,\mu_j}) \over p_{\eta_j}(x_{j,\mu_j})}\ \ \ 
\label{eq-fact1}
\eeq
The last step uses the fact that for $j=1,\ldots,n$, 
$x_{j,\nu_j} = x_{j-1*j} = x_{j-1,\mu_{j-1}}$.

From Steps (1) and (2a), we
see that $\Pi_0(\nu_j) = 1\,/\,(K_j\!+\!1)$ and $\Pi_1(\mu_j) = 
1\,/\,(K_j\!+\!1)$.   Using this, and again using 
$x_{j,\nu_j} = x_{j-1,\mu_{j-1}}$, we get that
\beq
   \lefteqn {{ 
     \prod\limits_{j=0}^{n-1} \Pi_1(\mu_j)\,
     \prod\limits_{j=1}^{n} \Pi_1(\nu_j\,|\,x_j)
     \over 
     \prod\limits_{j=1}^{n} \Pi_0(\nu_j)\,
     \prod\limits_{j=0}^{n-1} \Pi_0(\mu_j\,|\,x_j)}
   \ \ = \ \
   { \prod\limits_{j=0}^{n-1} \Pi_1(\nu_{j+1}\,|\,x_{j+1})\,(K_{j+1}\!+\!1)
     \over 
     \prod\limits_{j=0}^{n-1} \Pi_0(\mu_j\,|\,x_j)\,(K_j\!+\!1)
   }}\ \ \ \ \ \ \ \ \\[5pt]
   & = &
   \prod_{j=0}^{n-1}\,\
   {\displaystyle
    \ {p_{j*j+1}(x_{j+1,\nu_{j+1}}) \over p_{\eta_{j+1}}(x_{j+1,\nu_{j+1}})}
    \ \Big/\ {1 \over K_{j+1}\!+\!1}
    \sum_{k=0}^{K_{j+1}} {p_{j*j+1}(x_{j+1,k}) \over p_{\eta_{j+1}}(x_{j+1,k})}\
   \over\displaystyle
    {p_{j*j+1}(x_{j,\mu_j}) \over p_{\eta_j}(x_{j,\mu_j})}
    \ \Big/\ {1 \over K_j\!+\!1}
    \sum_{k=0}^{K_j} {p_{j*j+1}(x_{j,k}) \over p_{\eta_j}(x_{j,k})}
   } \\[5pt]
  & = &
   \prod_{j=0}^{n-1}\, 
    {p_{\eta_j}(x_{j,\mu_j}) \over p_{\eta_{j+1}}(x_{j,\mu_j})}\ 
   \prod_{j=0}^{n-1}\, 
   \left[
    {1 \over K_j\!+\!1}
    \sum_{k=0}^{K_j} {p_{j*j+1}(x_{j,k}) \over p_{\eta_j}(x_{j,k})}
    \ \Big/\
    {1 \over K_{j+1}\!+\!1}
    \sum_{k=0}^{K_{j+1}} {p_{j*j+1}(x_{j+1,k}) \over p_{\eta_{j+1}}(x_{j+1,k})} 
   \right]\ \ \ \ 
\label{eq-fact2}
\eeq

From Steps (1) and (3), we see that 
$\Pi_0(\nu_0) = \Pi_1(\nu_0) = 1\,/\,(K_0\!+\!1)$ and 
$\Pi_1(\mu_n) = \Pi_0(\mu_n) = 1\,/\,(K_n\!+\!1)$, so these factors 
cancel in equation~(\ref{eq-ratio01}).  The factors in
equation~(\ref{eq-fact1}) cancel with the first part of 
equation~(\ref{eq-fact2}).  The final result is that the simple importance 
sampling estimate based on a single LIS run is as shown in 
equation~(\ref{eq-lis}), demonstrating that $\rhatLIS$ is indeed an unbiased 
estimate of $r=Z_1/Z_0$.

\subsection{\hspace*{-4pt}Bridged LIS estimates}\vspace*{-4pt}

Since the LIS estimate can be viewed as a simple importance sampling
estimate on an extended space, we can consider a `bridged LIS'
estimate in which this top-level SIS estimate is replaced by a bridge
sampling estimate.  This will require that we actually perform the reverse
LIS procedure described above, from which an  LIS estimate for
the reverse ratio, $\underline{r} = Z_0/Z_1$, can be computed:
\beq
   \rhatLISrev^{(i)} & = & \prod_{j=1}^{n} \left[
     {1 \over K_j+1}\, \sum_{k=0}^{K_j}\,
          { p_{j-1*j}(x_{j,k}) \over p_{\eta_j}(x_{j,k}) }
     \ \Big/\
     {1 \over K_{j-1}+1}\, \sum_{k=0}^{K_{j-1}}\, 
          { p_{j-1*j}(x_{j-1,k}) \over p_{\eta_{j-1}}(x_{j-1,k}) }
     \right]
\label{eq-lis-rev}
\eeq
The reversed procedure requires independent sampling from $\pi_1$.  
This will usually not be possible directly, but well-separated states
from a Markov chain sampler with $\pi_1$ as its invariant distribution will 
provide a good approximation, provided that this sampler moves around the
whole distribution, without being trapped in an isolated mode.  Indeed,
the entire sample of $K_n\!+\!1$ states from $\pi_1$ that is needed
at the start of the reverse procedure can be obtained by taking consecutive 
states from such a Markov chain sampler.

For the bridged form of LIS, we also need a suitable bridge
distribution, $P_*$, for which we must be able to evaluate the ratios 
$P_*/P_0$ and $P_*/P_1$.  (Note that this choice of a
`top-level' bridge distribution is separate from the choices of
`low-level' bridge distributions, $p_{j*j+1}$, though we might use the same
form for both.)  With the optimal bridge of
equation~(\ref{eq-opt-bridge}), these ratios can be written as follows,
if the forward procedure is performed $M$ times and the reverse procedure
$\underline{M}$ times:
\beq
 {P\opt_*(x,\mu,\nu) \over P_0(x,\mu,\nu)} & = &
 \left[\,r\,(M/\underline{M})\,
       \left({P_1(x,\mu,\nu) \over P_0(x,\mu,\nu)}\right)^{-1}
       \!\! +\ 1\,\right]^{-1}
 \\[6pt]
 {P\opt_*(x,\mu,\nu) \over P_1(x,\mu,\nu)} & = &
 \left[\, r\,(M/\underline{M})\ +\ 
       \left({P_0(x,\mu,\nu) \over P_1(x,\mu,\nu)}\right)^{-1}
       \right]^{-1}
\eeq
The geometric bridge of equation~(\ref{eq-geo-bridge}) results in
\beq
 {P\geo_*(x,\mu,\nu) \over P_0(x,\mu,\nu)} & = &
 \sqrt{P_1(x,\mu,\nu) \over P_0(x,\mu,\nu)}
 \\[6pt]
 {P\geo_*(x,\mu,\nu) \over P_1(x,\mu,\nu)} & = &
 \sqrt{P_0(x,\mu,\nu) \over P_1(x,\mu,\nu)}
\eeq
These expressions allow us to express bridged LIS estimates in terms
of the simple LIS estimate of equation~(\ref{eq-lis}), and its reverse
version of equation~(\ref{eq-lis-rev}).  For the optimal bridge, we get
\beq
 \rhatLISbridged\opt & = & 
   {1 \over M} \sum_{i=1}^M\, 
      {1 \over r\,(M/\underline{M})\,/\,\rhatLIS^{(i)}\ +\ 1}
   \,\ \Big/\ 
   {1 \over \underline{M}} \sum_{i=1}^{\underline{M}}\, 
      {1 \over r\,(M/\underline{M})\ +\ 1/\rhatLISrev^{(i)}}
\label{eq-bridged-lis1}
\eeq
Similarly, for the geometric bridge, we get
\beq
 \rhatLISbridged\geo
   & = & {1 \over M} \sum_{i=1}^M\, \sqrt{\rhatLIS^{(i)}} \,\ \Big/\ 
         {1 \over \underline{M}} \sum_{i=1}^{\underline{M}}\, 
           \sqrt{\rhatLISrev^{(i)}}
\label{eq-bridged-lis2}
\eeq

\subsection{\hspace*{-4pt}LIS estimates with independent sampling with no
                          intermediate distributions}\vspace*{-4pt}

It is interesting to look at the special case of Linked Importance
Sampling with $n=1$ --- ie, in which the are no intermediate
distributions between $\pi_0$ and $\pi_1$ --- in which the points from both
$\pi_0$ and $\pi_1$ are sampled independently.  The LIS procedure
can then be simplified somewhat, and it is also possible to improve
the LIS estimate by averaging over the choice of link state.  Such
averaging is not feasible when Markov chain sampling is used, since
choosing a different link state would require a new simulation of the
Markov transitions.  

Since we will sample points independently, there is no need to decide
how many points will be sampled by the forward transitions and how
many by the reverse transitions in Steps (2a) and (2b) of the LIS
procedure.  We simply obtain a pair of samples consisting of
points $x_{0,0},\ldots,x_{0,K_0}$ drawn independently from $\pi_0$,
and points $x_{1,1},\ldots,x_{1,K_1}$ drawn independently from
$\pi_1$.  We then randomly select a link state, indexed by $\mu$, from 
among $x_{0,0},\ldots,x_{0,K_0}$ according to the 
following probabilities, which depend on the choice of a single
bridge distribution, denoted by $p_*(x)$:
\beq
  \Pi_0 (\mu \,|\, x_0) & = & 
  { p_*(x_{0,\mu}) \over p_0(x_{0,\mu})}\ \Big/\
               \sum\limits_{k=0}^{K_0} {p_*(x_{0,k}) \over p_0(x_{0,k}) }
\eeq
The LIS estimate for $r = Z_1/Z_0$ based on this pair of samples
from $\pi_0$ and $\pi_1$ is
\beq
 \rhatLIS^{(i)} & = &
 {1 \over K_0\!+\!1} \sum_{k=0}^{K_0} {p_*(x_{0,k}) \over p_0(x_{0,k})}
 \ \Big/\,
 {1 \over K_1\!+\!1} \left[ {p_*(x_{0,\mu}) \over p_1(x_{0,\mu})} 
         \, +\, \sum_{k=1}^{K_1} {p_*(x_{1,k}) \over p_1(x_{1,k})} 
 \right]
\label{eq-lis-indep}
\eeq
The superscript $i$ is used here
to indicate that this estimate is based on the $i$'th 
pair of samples.  We can see that it is very similar to the bridge sampling 
estimate of equation~(\ref{eq-bridge}), except that the link state is included 
in both samples.  Since these LIS estimates are unbiased, we can 
average $M$ of them to obtain a final LIS estimate.

We can also average the estimate of equation~(\ref{eq-lis-indep})
over the random choice of link state, which
is guaranteed to produce an estimate (also unbiased) with smaller 
mean-squared-error (see Schervish 1995, Section 3.2).  The result is
\beq
 \rhatLISave^{(i)} & = & 
 \sum_{\mu=0}^{K_0} \Pi_0(\mu\,|\,x_0) \
 {1 \over K_0\!+\!1} \sum_{k=0}^{K_0} {p_*(x_{0,k}) \over p_0(x_{0,k})}
 \ \Big/\,
 {1 \over K_1\!+\!1} \left[ {p_*(x_{0,\mu}) \over p_1(x_{0,\mu})}
         \, +\, \sum_{k=1}^{K_1} {p_*(x_{1,k}) \over p_1(x_{1,k})} 
 \right] \\[5pt]
& = &
 {K_1\!+\!1 \over K_0\!+\!1}\ \sum_{\mu=0}^{K_0} \
 {p_*(x_{0,\mu}) \over p_0(x_{0,\mu})}
 \ \Big/\,
  \left[ {p_*(x_{0,\mu}) \over p_1(x_{0,\mu})}
         \, +\, \sum_{k=1}^{K_1} {p_*(x_{1,k}) \over p_1(x_{1,k})} 
 \right]
\label{eq-lis-ave}
\eeq
Averaging these estimates over $M$ pairs of samples produces a final estimate
denoted by $\rhatLISave$.

To use bridged LIS in this context, we need to find reverse estimates
as well, but these reverse estimates needn't be independent of the
forward estimates, since the asymptotic validity of the bridge
sampling estimate of equation~(\ref{eq-bridge}) does not depend on the
samples $x_0$ and $x_1$ being independent.  Accordingly, we can use
the same samples from $\pi_0$ and $\pi_1$ for the forward and the
reverse operations.  However, to perform reverse sampling, we need to
have a sample of $K_1\!+\!1$ points drawn from $\pi_1$, the first of
which is ignored when performing forward sampling.  Conversely, the
first of the $K_0\!+\!1$ points drawn from $\pi_0$ is ignored when
performing the reverse sampling.

We can improve the bridged LIS estimates by averaging the numerator
and the denominator of equation~(\ref{eq-bridged-lis1})
or~(\ref{eq-bridged-lis2}) with respect to the random choice of link
state.  We can also average with respect to the omission of one of the
points from one of the samples --- ie, rather than omitting the first
of $K_1 + 1$ points in the sample from $\pi_1$ when computing a
forward estimate, we average with respect to a random choice of point
to omit, and similarly for reverse estimates.  Note that the averaging
should be done over the sums in the numerator and denominator, not
with respect to the entire estimate, nor with respect to the values of
$\rhatLIS^{(i)}$ and $\rhatLISrev^{(i)}$ appearing inside the
summands.  The effective sample size after this additional averaging
of dependent points is unclear, so it is not obvious what the ratio
of sample sizes in equation~(\ref{eq-opt-bridge}) should be, but
using $(K_0\!+\!1)/(K_1\!+\!1)$ is probably adequate.

\section{\hspace*{-7pt}Analytical comparisons of AIS and 
                       LIS}\label{sec-anal}\vspace*{-10pt}

In this section, I analyse (somewhat informally) the performance of
AIS and LIS asymptotically, and in other situations where analytical
results are possible.

\subsection{\hspace*{-4pt}Asymptotic properties of 
                       AIS and LIS estimates}\label{sec-asym}\vspace*{-4pt}

I begin by analysing the asymptotic performance of AIS and LIS when
the sequence of distributions is defined by an unnormalized density function
of the following form:
\beq
   p_{\eta}(x) & = & p_0(x)\, \exp (-\eta U(x))
\label{eq-U-dist}
\eeq
This class includes sequences of canonical distributions defined by 
equation~(\ref{eq-canonical}) in which the inverse temperature
varies, as well as
sequences that can be used for Bayesian analysis, in which $p_0$ defines the
prior and $\eta$ is a power that the likelihood (expressed as $\exp(-U(x))$) is 
raised to, with $\eta=1$ giving the posterior distribution.  
For these distributions, we can express $r$ using the well-known 
`thermodynamic integration' formula as follows:
\beq
  r\ \ =\ \ \log(Z_1/Z_0)\ \ =\ \ - \int_0^1 E_{\pi_{\eta}}(U)\,d\eta
\label{eq-therm-int}
\eeq

The analysis here is asymptotic, as the number of intermediate
distributions used, given by $n\!-\!1$, goes to infinity.  I will
assume the $\eta_j$ defining these distributions are chosen according to a 
scheme in which for any 
$a \in (0,1)$, the spacing $\eta_{j+1}-\eta_j$ when $j = \lfloor a\,n \rfloor$ 
is asymptotically proportional to $1/n$ --- in other words, 
the relative density of intermediate distributions in the neighborhood
of different values of $\eta$ stays the same as the overall density increases. 
The simplest such scheme is to let $\eta_j = j/n$, though other schemes 
may sometimes be better.

With the above form for $p_{\eta}$, the AIS estimate from a single run 
(from equation~(\ref{eq-ais-est})) can be written as follows:
\beq
  \log\ \rhatAIS^{(i)}
  & = & 
   \sum_{j=0}^{n-1}\, \log \Big(p_{\eta_{j+1}}(x^{(i)}_j)
                                \,\Big/\,p_{\eta_j}(x^{(i)}_j)\Big)
   \ \ =\ \  
   \sum_{j=0}^{n-1}\, - (\eta_{j+1}-\eta_j)\, U \Big(x^{(i)}_j\Big)
\label{eq-ais-reim}
\eeq
When $\eta_j=j/n$, this can be seen as a stochastic form of Riemann's Rule 
for numerically integrating equation~(\ref{eq-therm-int}), though one 
difference is that $\log\ \rhatAIS$ converges to the correct value as $M$ goes 
to infinity even if $n$ stays fixed.

Provided that there is some finite bound on the variance of $U$ under all 
the distributions $\pi_{\eta}$, and that the Markov transitions used mix well,
a Central Limit Theorem will apply, allowing us to conclude that the
distribution of $\ell_n = \log\ \rhatAIS^{(i)}$ becomes
Gaussian as $n$ goes to infinity.  Let the mean of $\ell_n$ be $\mu_n$,
and let the variance of $\ell_n$ asymptotically be $\sigma^2/n$, where $\sigma$
is determined by details of the spacing of intermediate distributions and
of the degree of autocorrelation in the Markov transitions.
Note that $E[Y^q]=\exp(q\mu+q^2\varsigma^2/2)$ when $Y=\exp(X)$
and $X$ is Gaussian with mean $\mu$ and variance $\varsigma^2$.  
Using this, the mean of $\exp(\ell_n)$ is $\exp(\mu_n+\sigma^2/2n)$.  This
must equal $r$, since $\rhatAIS$ is unbiased, so $\mu_n = \log(r)-\sigma^2/2n$.
Using this, we can see that the variance of $\rhatAIS^{(i)}=\exp(\ell_n)$ is 
$r\,[\exp(\sigma^2/2n) - 1]$, which for large $n$ will be approximately
$r\sigma^2/2n$.  The variance of $\rhatAIS$ will therefore be $r\sigma^2/2nM$.
Asymptotically, the total computational effort, which will generally be
proportional to $nM$, can be divided in any way between more intermediate
distributions ($n$) or more runs ($M$) without affecting the accuracy
of estimation of $r$, provided that $n$ is kept large enough that
these asymptotic results apply --- a fact noted by Hendrix and Jarzynski (2001).
We can therefore use a value of $M$ greater than one without penalty, 
in order to obtain an error estimate from the degree of variation 
over the $M$ runs.

For LIS, we can write the log of the estimate from one run 
(equation~(\ref{eq-lis})) as follows:
\beq
   \log\ \rhatLIS^{(i)} & = & \sum_{j=0}^{n-1} \left[
     \log \left({1 \over K_j+1}\, \sum_{k=0}^{K_j}\,
          { p_{j*j+1}(x_{j,k}) \over p_{\eta_j}(x_{j,k}) } \right)
     \ -\
     \log \left({1 \over K_{j+1}+1}\, \sum_{k=0}^{K_{j+1}}\, 
          { p_{j*j+1}(x_{j+1,k}) \over p_{\eta_{j+1}}(x_{j+1,k}) }\right)
     \right]\ \ \ \ \ \
\label{eq-logrhatLIS}
\eeq
Suppose that we let $K_j = \lceil m K_j^0 \rceil$ for all $j$ and some set of
$K^0_j$, and that we then let $m$ go to infinity.  Assuming that the variances
of the ratios of probabilities are finite, and that the Markov chain transitions
used mix sufficiently well, a Central Limit
Theorem will again apply, and we can conclude that all of the $n$ terms in 
the sum above, and therefore also the sum itself, will approach Gaussian 
distributions, with variances proportional to $1/m$.  

To analyse the LIS estimate in more detail, we need to assume a form of
bridge distribution, as well as a form for $p_{\eta}$.  If $p_{\eta}$
has the form of equation~(\ref{eq-U-dist}) and we use the geometric bridge
of equation~(\ref{eq-geo-bridge}), we can write 
\beq
   \log\ \rhatLIS^{(i)} & = & \sum_{j=0}^{n-1}\, \left[\
     \log \left( {1 \over K_j+1}\, \sum\limits_{k=0}^{K_j}\,
          \exp(-(\eta_{j+1}\!-\!\eta_j)\, U(x_{j,k})\, /\, 2)  \right) 
   \ -\ \right. \nonumber \\[4pt]
   & & \ \ \ \ \ \ \ \ \left.
     \log \left( {1 \over K_{j+1}+1}\,\sum\limits_{k=0}^{K_j}\,
          \exp(-(\eta_j\!-\!\eta_{j+1})\, U(x_{j+1,k})\, /\, 2) \right)
     \ \right]
\eeq
Since $\exp(z)\approx1+z$ and $\log(1+z)\approx z$ when $z$ is small, we can
rewrite this when $n$ is large (and hence $\eta_{j+1}\!-\!\eta_j$ is small) as
\beq
   \log\ \rhatLIS^{(i)} & \approx & \sum_{j=0}^{n-1}\, \left[\
     \log \left( 1 \ -\ { \eta_{j+1}\!-\!\eta_j \over 2}\,
         {1 \over K_j+1}\, \sum\limits_{k=0}^{K_j} U(x_{j,k}) \right) 
   \ -\ \right. \nonumber \\[4pt]
   & & \ \ \ \ \ \ \ \ \left.
     \log \left( 1 \ +\ {\eta_{j+1}\!-\!\eta_j \over 2}\,
         {1 \over K_{j+1}+1}\, \sum\limits_{k=0}^{K_{j+1}} U(x_{j+1,k}) \right) 
     \ \right] \\[6pt]
   & \approx & \sum_{j=0}^{n-1}\,
         - {\eta_{j+1}\!-\!\eta_j \over 2}\,
         \left[ {1 \over K_j+1}\, \sum\limits_{k=0}^{K_j} U(x_{j,k})
         \ +\   {1 \over K_{j+1}+1}\, 
                \sum\limits_{k=0}^{K_{j+1}} U(x_{j+1,k}) \right]  \\[6pt]
   & = & 
         -\ {\eta_1\!-\!\eta_0 \over 2}\,
           {1 \over K_0+1}\, \sum\limits_{k=0}^{K_0} U(x_{0,k})
         \ -\ {\eta_n\!-\!\eta_{n-1} \over 2}\,
              {1 \over K_n+1}\, \sum\limits_{k=0}^{K_n} U(x_{n,k}) 
         \nonumber \\[4pt]
   &   &  -\ \sum_{j=1}^{n-1}\, 
           {\eta_{j+1}\!-\!\eta_{j-1} \over 2}\,
           {1 \over K_j+1}\, \sum\limits_{k=0}^{K_j} U(x_{j,k})
\eeq

When $\eta_j=j/n$, this looks like a stochastic form of the
Trapezoidal Rule for numerically integrating
equation~(\ref{eq-therm-int}).  Since the Trapezoidal Rule converges
faster than Reimann's Rule, one might expect LIS to perform better
than AIS asymptotically, but this is not so in this stochastic
situation.  Suppose for simplicity that we set all $K_j=m$.  The
variance of $\log\ \rhatLIS^{(i)}$ will be dominated by the variance
of the last sum above, which will be proportional to $1/nm$, assuming
that $m$ is large, so that the dependence between terms (from sharing
link states) is negligible.  Using the same argument as for AIS above,
the variance of $\log \rhatLIS$ will be proportional to $1/nmM$.
Considering that the computation time for an LIS run will be
proportional to $nm$, versus $n$ for AIS, we see that the variances of
the AIS and LIS estimates go down the same way in proportion to
computation time, asymptotically as $n$ and $m$ go to infinity.

Furthermore, the proportionality constant should be the same for
AIS and LIS, assuming that the overhead of the two procedures is
negligible compared to the time spent performing Markov transitions,
so that the proportionality constants for computation time are the
same for AIS (multiplying $n$) and for LIS (multiplying $nm$).  The
proportionality constants for variance for AIS (multiplying $1/nM$)
and for LIS (multiplying $1/nmM$) depend in a complex way on the form of
the density of $\eta_j$ values and on the mixing properties of the
Markov transitions, but the result should be the same for AIS and
LIS, provided the same scheme is used for choosing $\eta_j$ values,
and the same Markov transitions are used, parameterized smoothly in
terms of $\eta$.  A difference that might appear significant is that
for AIS only one Markov transition is done for each $\eta_j$, whereas
for LIS, $m$ such transitions are done.  However, as $n$ goes to
infinity, nearby distributions become more similar, so transitions for
$m$ consecutive distributions become similar to $m$ transitions for
one of these distributions.

The apparently pessimistic conclusion from this is that when both $n$
and $m$ (and hence the $K_j$) are large, the performance of LIS should
be about the same as that of AIS (with $n$ for AIS chosen to equalize
the computation time), assuming that the distributions used have the
form of equation~(\ref{eq-U-dist}), that the variance of $U$ is finite
under all of the distributions $\pi_{\eta}$, and that the Markov
transitions used mix well enough.  Fortunately, however, there is no
reason to make both $m$ and $n$ large with LIS.  For good performance,
$n$ must be large enough that $\pi_{\eta_j}$ and $\pi_{\eta_{j+1}}$
overlap significantly, but there is no reason to make $n$ much larger
than this.  The accuracy of the estimates can be improved as desired
by increasing $m$ and/or $M$ while keeping $n$ fixed.  The results
below show that LIS estimates with $n$ fixed are sometimes much better
than AIS estimates.

Finally, let us consider the asymptotic performance of the bridged
versions of AIS and LIS, assuming that the variance of $U$ is finite,
so that the distribution of the estimates from individual runs becomes
Gaussian as $n$ (for AIS) or $m$ (for LIS) goes to infinity.  Looking
at equations~(\ref{eq-bridged-lis1}) and~(\ref{eq-bridged-lis2}),
which also are applicable to bridged AIS estimates, we see that the
log of $\rhatLISbridged^{(i)}$ can for both optimal and geometric
bridges be expressed as the difference of the log of the numerator,
which is the mean of a function of the forward estimates,
$\rhatLIS^{(i)}$, and the log of the denominator, which is the mean of
a function of the reverse estimates, $\rhatLISrev^{(i)}$.  If these
forward and reverse estimates have Gaussian distributions with small
variances, $\sigma^2$ and $\underline{\sigma}^2$, then
$\rhatLISbridged^{(i)}$ will also be Gaussian, with a variance that
can be computed in terms of the derivatives of the summands in the
numerator and the denominator, with respect to $\rhatLIS^{(i)}$ and
$\rhatLISrev^{(i)}$, evaluated at the true values of $r$ and $1/r$.
I will assume that $r=1$ below, as can be done without loss of generality.

For the geometric bridge, these derivatives are both $1/2$, from which
it follows that the variance of the numerator in
equation~(\ref{eq-bridged-lis2}) is $\sigma^2/4M$ and that of the
denominator is $\underline{\sigma}^2/4\underline{M}$.  Since the
numerator and denominator evaluate to one for $\rhatLIS^{(i)}=r=1$ and
$\rhatLISrev^{(i)}=1/r=1$, the sum of the variances of the logs of the
numerator and denominator is $\sigma^2/4M +
\underline{\sigma}^2/4\underline{M}$. If
$\sigma^2=\underline{\sigma}^2$ and $M=\underline{M}$, this reduces to
$\sigma^2/2M$.  The variance of an unbridged LIS estimate will be
$\sigma^2/M$.  However, the bridged estimate requires time
proportional to $M+\underline{M}$, compared to just $M$ for the
unbridged estimate.  The value of $M$ for the unbridged method can
therefore be twice as large as for the bridged method, with the result
that bridged and unbridged estimates perform equally well
asymptotically (assuming the variance of $U$ is finite).

For the optimal bridge, the derivatives of the summands in the
numerator and denominator are both $1/4$, when evaluated at
$\rhatLIS^{(i)}=r=1$ and $\rhatLIS^{(i)}=1/r=1$, and assuming that
$M=\underline{M}$.  The numerator and denominator both evaluate to
$1/2$, with the result that asymptotically the variance of the bridged
estimate, assuming $\sigma^2=\underline{\sigma}^2$, is $\sigma^2/2M$,
the same as for the geometric bridge.

In conclusion, bridged AIS and LIS estimates asymptotically have the
same performance as the corresponding unbridged estimates (with twice
the value of $M$), for both the optimal and geometric bridges,
assuming $U$ has finite variance.  This conclusion applies more
generally, as long as a Central Limit Theorem holds for the individual
estimates, $\rhatLIS^{(i)}$ and $\rhatLISrev^{(i)}$.  However, the
bridged methods may be much better when the variance of $U$ is
infinite, or for classes of distributions other than that of
equation~(\ref{eq-U-dist}).  The bridged methods may also provide
improvement when the values of $n$ or $m$ are not large enough for the
asymptotic results to apply.

\subsection{\hspace*{-4pt}Properties of AIS and LIS when sampling from
            uniform distributions}\label{sec-unif}\vspace*{-4pt}

In this section, I will demonstrate that when $n$ is kept suitably
small, LIS can perform much better than AIS when these methods are
applied to sequences of uniform distributions.

As a first example, consider the class of nested uniform 
distributions with unnormalized densities given by\vspace*{-6pt}
\beq
  p_{\eta}(x) & = & \left\{ \begin{array}{ll} 
      1 & \mbox{if $-s^{\eta} < x < s^{\eta}$} \\ 0 & \mbox{otherwise} 
  \end{array}\right.
\eeq
for which the normalizing constants are $Z_{\eta} = 2s^{\eta}$, so that
$r = Z_1/Z_0 = s$.  The results concerning this class of distributions
can easily be extended to any class of uniform distributions, in any
number of dimensions, that have nested regions of support.
For both AIS and LIS, I will assume that the intermediate
distributions are defined by $\eta_j = j/n$.  With this choice, the
probability that a point, $x$, randomly sampled from $\pi_j$ will have
$p_{j+1}(x)=1$ is $s^{1/n}$, for any $j$.

During an AIS run, only a single point is sampled from each
distribution.  An AIS run will produce an estimate for $r$ of zero if
any of the ratios ${p_{\eta_{j+1}}(x^{(i)}_j)\,/\,
p_{\eta_j}(x^{(i)}_j)}$ in equation~(\ref{eq-ais-est}) are zero, which
happens with probability $1 - (s^{1/n})^n\, =\, 1-s$, and will
otherwise produce an estimate of one.  Note that the distribution of
estimates is independent of $n$.  AIS is therefore not a useful
technique for nested uniform distributions --- simple importance
sampling (ie, AIS with $n\!=\!1$) would work just as well (or just as
poorly, if $s$ is very small).  Bridged AIS produces no improvement in
this context.

Suppose instead we use LIS with all $K_j=m$, and suppose that the
Markov transitions, $T_j$, produce points that are almost independent
of the previous point.  For this problem, both the geometric and
optimal forms of the bridge distribution result in $p_{j*j+1}(x) =
p_{\eta_{j+1}}(x)$.  If $m+1$ points are sampled independently from
$\pi_{\eta_j}$, the fraction of these points for which
$p_{\eta_{j+1}}(x)$ is one will have variance
$s^{1/n}\,(1\!-\!s^{1/n})\,/\,(m\!+\!1)$.  For sufficiently large
$m$, the variance of the log of this fraction will be
approximately $(s^{1/n}\,(1\!-\!s^{1/n})\,/\,(m\!+\!1))\,/\,s^{2/n}$,
which simplifies to $(s^{-1/n}\!-\!1)\,/\,(m\!+\!1)$.  For this
approximation to be useful, the probability that none of the $m+1$
points sampled from $\pi_{\eta_j}$ lie in the region where
$p_{\eta_{j+1}}$ is one, equal to $(1-s^{1/n})^{m+1}$, must be negligible.  
This probability must be fairly small anyway, if LIS is to perform well.

Suppose that the computational cost of an LIS run is proportional to
the sum of the number of points sampled from $\pi_0$ and the number of
Markov transitions performed.  If we fix this cost, the number of
intermediate distributions, $n$, and the number of transitions for
each distribution, $m$, will be related by $m(n\!+\!1)\,=\,C$, for
some constant $C$.  Assume for the moment that both $n$ and $m$ are
large.  The probability of a run producing a zero estimate will then
be negligible, and we can assess the accuracy of the estimate for one
run by the variance of $\log \rhatLIS^{(i)}$ (modified in some way
to eliminate the infinity resulting from the negligible, but non-zero,
probability that $\rhatLIS^{(i)}$ is zero).  Looking at
equation~(\ref{eq-logrhatLIS}), we see that for these nested uniform
distributions, the second log term vanishes ---
$p_{j*j+1}(x_{j+1,k})\,/\,p_{\eta_{j+1}}(x_{j+1,k})$ is always one,
since $p_{j*j+1}$ is the same as $p_{\eta_{j+1}}$.  When $m$ is large,
the dependence between terms with different values of $j$ will be
negligible, so we can add the variances of the terms to get the variance
of the estimate, obtaining the result that
\beq
  \Var \Big(\log\ \rhatLIS^{(i)}\Big) 
    & \approx & n\,(s^{-1/n}\!-\!1)\,/\,(m\!+\!1)
\label{eq-varLIS-nest}
\eeq
When $n$ is large, $s^{-1/n}=\, \exp(\log(1/s)/n)$ is approximately
$1+\log(1/s)/n$, and hence the variance above is
approximately $\log(1/s)\,/\,(m\!+\!1)$.
So it seems that the larger the value of $m$, the better ---
until we reach a value of $m$ for which the corresponding value of $n$,
equal to $C/m\,-\,1$, is small enough that this result no longer applies.

Best performance will therefore come using a fairly small value of
$n$, but a large value of $m$.  Substituting $m=C/(n\!+\!1)$ into
equation~(\ref{eq-varLIS-nest}), and assuming $m/(m\!+\!1)\approx 1$, we get
\beq
  \Var \Big(\log\ \rhatLIS^{(i)}\Big) 
    & \approx & n\,(s^{-1/n}\!-\!1)\,/\,(C/(n\!+\!1))
    \ \ =\ \ n(n\!+\!1)\,(s^{-1/n}\!-\!1)\,/\,C
\eeq
The value of $n$ that minimizes this depends only on $s$, not on $C$.
The optimal choice of $n$ increases slowly as $s$ gets smaller:\ \
$s=0.1$ gives $n=2$, $s=0.05$ gives $n=3$, $s=0.01$ gives $n=4$, and
$s=0.0001$ gives $n=7$.

As a second example, consider the class of non-nested uniform distributions 
with unnormalized densities given by\vspace*{-6pt}
\beq
  p_{\eta}(x) & = & \left\{ \begin{array}{ll} 
      1 & \mbox{if $\eta t-1 < x < \eta t+1$} \\ 0 & \mbox{otherwise} 
  \end{array}\right.
\eeq
For this class, $Z_{\eta} = 2$ for all $\eta$, so $r = Z_1/Z_0 = 1$.
I will again assume that the intermediate
distributions are defined by $\eta_j = j/n$, and that all $K_j=m$.  Assuming 
that $n$ is greater than
$t/2$, the probability that a point, $x$, randomly sampled from $\pi_{\eta_j}$ 
will have $p_{\eta_{j+1}}(x)=1$ is $1-t/2n$, for any $j$.

For this example, AIS estimates do not converge to the true value of
$r$ as $M$ increases, regardless of the value of $n$.  To see this,
note that the ratios in equation~(\ref{eq-ais-est}) will all be either
zero or one, and that the estimate from one run, $\rhatAIS^{(i)}$,
will be one if all of these ratios are one, and zero otherwise.  The
probability of a particular ratio being one is $1-t/2n$, so the
probability that all are one (assuming the $T_{\eta}$ produce points
independent of the current point) is $(1-t/2n)^n$, which approaches
$\exp(-t/2)$ as $n$ goes to infinity.  The AIS estimate, averaging
over $M$ runs, will have mean $\exp(-t/2)$, rather than the correct
value of one.

In contrast, bridged AIS estimates will converge to the true value as $M$
increases, as long as $n$ is at least $t/2$, so that there is overlap
between successive distributions in the sequence.  However, when $t$
is large, the overlap between the distributions over paths produced by
forward and reverse AIS runs, given by $\exp(-t/2)$, will be very
small, and the procedure will be very inefficient.

To see how well LIS performs, recall the formula for $\log \rhatLIS$
from equation~(\ref{eq-logrhatLIS}):  
\beq
   \log\ \rhatLIS^{(i)} & = & \sum_{j=0}^{n-1} \left[
     \log \left({1 \over K_j+1}\, \sum_{k=0}^{K_j}\,
          { p_{j*j+1}(x_{j,k}) \over p_{\eta_j}(x_{j,k}) } \right)
     \ -\
     \log \left({1 \over K_{j+1}+1}\, \sum_{k=0}^{K_{j+1}}\, 
          { p_{j*j+1}(x_{j+1,k}) \over p_{\eta_{j+1}}(x_{j+1,k}) }\right)
     \right]\ \ \ \ \ \
\eeq
Due to symmetry, the two log terms above have the same distribution,
for all $j$.  The variance of one of these log
terms (for large $m$) is
$((t/2n)\,(1\!-\!t/2n)\,/\,(m\!+\!1))\,/\,(1\!-\!t/2n)^2$, which
simplifies to $1\,/\,((2n/t\!-\!1)\,(m\!+\!1))$.  The second log
term in equation~(\ref{eq-logrhatLIS}) for one $j$ will involve the
same points, $x_{j+1,k}$, as the first log term for the next $j$.  The
effect of this is that these terms will be negatively correlated, with 
correlation of $-1$ if $n\!=\!t$.  However, since the
two terms occur with opposite signs, the effect on the final sum is 
that $n\!-\!1$ pairs of terms (out of $2n$ terms total) are positively
correlated.  Straightforward calculations show that this correlation is
$2n/t - 1$ for $t/2 < n \le t$ and $1\,/\,(2n/t - 1)$ for $n \ge t$.  
Using the fact that when $X$ and $Y$ have the same
distribution, $\Var(X+Y) = 2\,\Var(X)\,[1+\Cor(X,Y)]$, we obtain the 
result that, for large $m$,
\beq
  \Var \Big(\log\ \rhatLIS^{(i)}\Big) 
    & \approx & {2 \over (2n/t\!-\!1)\,(m\!+\!1)} 
                \left\{\begin{array}{ll}
                   n\ +\ (n\!-\!1)\,(2n/t-1)
                     & \ \ \mbox{if $t/2 < n \le t$}
                \\[4pt]
                   n\ +\ (n\!-\!1)\,/\, (2n/t-1)
                     & \ \ \mbox{if $n \ge t$} 
                \end{array}\right\} 
\eeq
Setting $m = C/(n\!+\!1)$, and assuming $m/(m\!+\!1)\approx 1$, gives
\beq
  \Var \Big(\log\ \rhatLIS^{(i)}\Big) 
    & \approx & {2(n\!+\!1) \over C(2n/t\!-\!1)} 
                \left\{\begin{array}{ll}
                   n\ +\ (n\!-\!1)\,(2n/t-1)
                     & \ \ \mbox{if $t/2 < n \le t$}
                \\[4pt]
                   n\ +\ (n\!-\!1)\,/\, (2n/t-1)
                     & \ \ \mbox{if $n \ge t$} 
                \end{array}\right\}
\eeq
Numerical investigation shows that the global minimum of the variance
occurs where $n$ is near $(3/2)\,t$.  A second local minimum where $n$
is near $(3/4)\,t$ also exists.  The two minima are nearly equally good
when $t$ is large.  There is a local maximum where $n$ is near $t$,
with the variance there being about 19\% greater than at the global
minimum.  The variance is much larger for very large and very small values
of $n$.  We therefore see that for this example too, the best results
are obtained by fixing $n$ to a moderate value; any desired level of
accuracy can then be obtained by increasing $m$ and/or $M$.

\section{\hspace*{-7pt}Empirical comparisons of AIS and 
                       LIS}\label{sec-cmp}\vspace*{-10pt}

The analytical results of the previous section indicate that LIS can
sometimes perform much better than AIS, but that the benefits of LIS
may only be seen when the number of intermediate distributions used is
kept suitably small (but not so small that they do not overlap).  In
this section, I investigate the performance of AIS and LIS (and their
bridged versions) empirically.  The programs used for these tests
(written in R) are available from my web page.

These tests were done using sequences of one-dimensional distributions 
having unnormalized density functions of the following form:
\beq
  p_{\eta}(x) & = & 
    \exp\Big(\!-\!\Big|(x\!-\!\eta t)\,/\,s^{\eta}\,\Big|^q\,\Big)
\eeq
where $s$, $t$, and $q$ are fixed constants.  As $\eta$ moves from 0 to 1,
the centre of this distribution shifts by $t$, and changes width by the
factor $s$.  The power $q$ controls how thick the tails of the distributions
are.  When $q=2$, the distributions are Gaussian; a larger value 
produces lighter tails.  Note that $Z_{\eta}$ is
proportional to $s^{\eta}$, and hence $r = Z_1/Z_0$ is equal to $s$.

If $t=0$, the distributions can be written in the form of
equation~(\ref{eq-U-dist}), after reparameterizing in terms of $\eta'
= 1/s^{\eta q}$, so that $p_{\eta'}(x) = \exp(-\eta' |x|^q)$.  In this
case, we expect the asymptotic behaviour to be as discussed in
Section~\ref{sec-asym}, but the behaviour with samples of practical
size may be different.  As $q$ goes to infinity, the distributions
converge to uniform distributions over $(\eta t\!-\!s^{\eta},\,\eta
t\!+\!s^{\eta})$, and the results of Section~\ref{sec-unif} become relevant.

\begin{figure}[t]

\vspace*{-29pt}

\centerline{\includegraphics{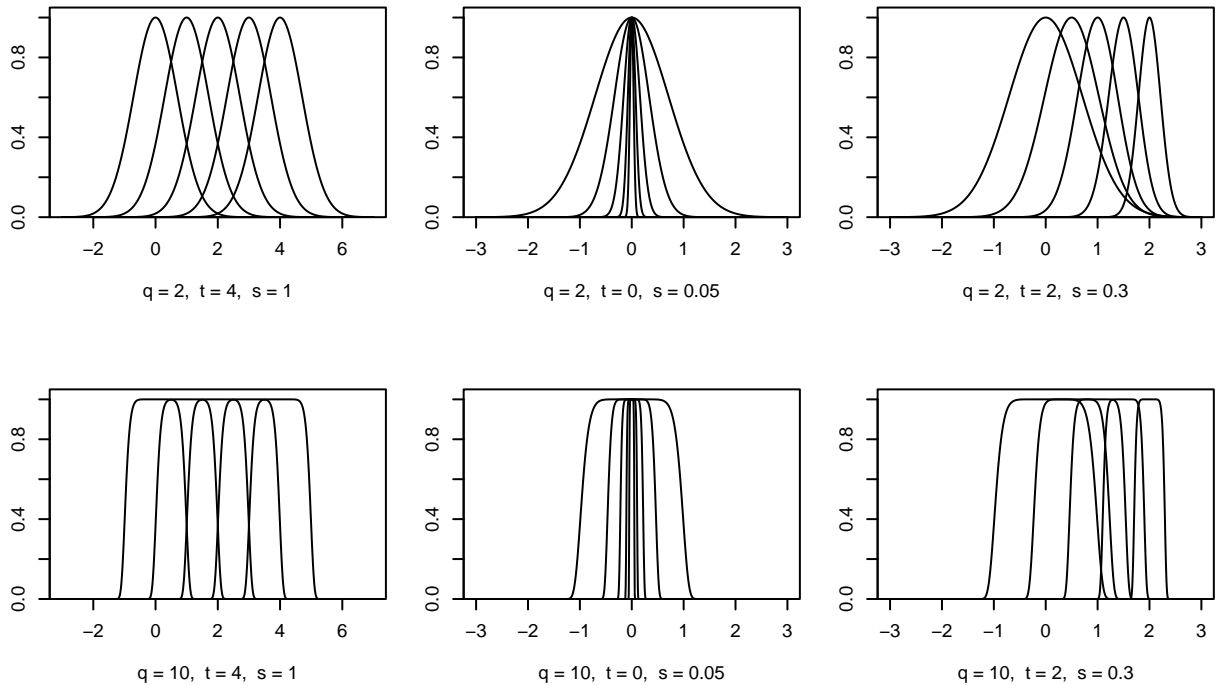}}

\caption[]{The sequences of unnormalized density functions used for the
           tests.  The plots show the unnormalized density functions for
           $\eta\, =\, 0,\, 1/4,\, 2/4,\, 3/4,\, 1$, for six combinations 
           of $s$, $t$, and $q$.}\label{fig-seq}

\end{figure}

I did an initial set of tests using six sequences of distributions.
Three of these sequences were of Gaussian distributions, with $q\!=\!2$.
The first of these used $s\!=\!1$ and $t\!=\!4$, producing a shift with no
change in scale as $\eta$ increases from 0 to 1.  The second used
$s\!=\!0.05$ and $t\!=\!0$, producing a contraction with no shift.  The last
used $s\!=\!0.3$ and $t\!=\!2$, combining a shift with a contraction.  A
second set of three sequences used the same values of $s$ and $t$, but
with $q\!=\!10$, which produces more `rectangular' distributions with
lighter tails.  The six sequences are shown in Figure~\ref{fig-seq}.
Each sequence in these plots consists of five distributions,
corresponding to $\eta\, =\, 0,\, 1/4,\, 2/4,\, 3/4,\, 1$.  These were
the sequences used for the LIS runs (hence $n\!=\!4$ for these runs).  The
AIS runs used more distributions, spaced more finely with respect to
$\eta$, so as to produce the same number of Markov transitions and
sampling operations as in the LIS runs.

These distributions (for any $\eta$) can easily be sampled from using
rejection sampling.  Samples from $\pi_0$ and $\pi_1$ were used to
initialize forward and reverse runs of AIS and LIS.  For this test, we
pretend that sampling for other $\pi_{\eta}$ must be done using Markov
chain methods.  The transition used for $\pi_{\eta}$, $T_{\eta}$, was
a random-walk Metropolis update, using a Gaussian proposal
distribution with mean equal to the current point and standard
deviation $s^{\eta}$.  Since Metropolis updates are reversible,
$\underline{T}_{\eta}$ was the same.

Two sets of forward and reverse LIS runs were done with $n\!=\!4$, all
$K_j\!=\!50$, and $M\!=\!20$, one set using the geometric bridge, the
other using the optimal bridge with the true value of $r$.  The
forward estimates were computed from equation~(\ref{eq-lis}); the
reverse estimates from equation~(\ref{eq-lis-rev}), which is
equivalent to using the forward procedure with the reverse sequence of
distributions.  Bridged LIS estimates were also found using
equation~(\ref{eq-bridged-lis1}), with the value of $r$ found by
iteration.  To make the comparison with forward and reverse estimates
fair, the bridged LIS estimates used $M\!=\!10$ --- ie, only half of
the forward and half of the reverse runs were used, for a total of
$20$ runs.

A corresponding set of forward, reverse, and bridged AIS runs were
also done, with $n\!=\!250$ and $M\!=\!20$ ($M\!=\!10$ for the bridged
estimates).  If sampling a point from $\pi_0$ or $\pi_1$ takes about
the same computation time as a Metropolis update, these AIS runs will
take about the same time as the LIS runs.  (This assumes that sampling
and Markov transitions dominate the time, which is typically true for
real problems but perhaps not for this simple test problem.)

Sets of longer LIS and AIS runs were also done, which were the same as
the sets above except that for LIS, $K_j\!=\!200$ for all $j$, and for
AIS, $n\!=\!1000$, which again equalizes the computation time.

Experience, together with the asymptotic results of
Section~\ref{sec-asym}, shows that estimates produced using a small
value of $M$ are better than, or at least as good as, those produced
with larger $M$.  I chose $M\!=\!20$ ($M\!=\!10$ for bridged estimates) since
this is about the smallest value that allows reliable estimation of
standard errors, which would usually be needed in practice.

The standard errors for AIS and LIS estimates of $\rhat$ were
estimated by the sample standard deviation of the $\rhat^{(i)}$
divided by $\sqrt{M}$.  When comparing the methods, I looked primarily
at the mean squared error when estimating $\log(r)$ (rather than when
estimating $r$).  The estimate I used was $\log(\rhat)$, and the
standard error for this estimate was estimated by the standard error
for $\rhat$ divided by $\rhat$.  For the reverse runs, $\log(r)$ was
estimated by $-\log(\rhatrev)$.  For bridged AIS and LIS, the standard
errors for the log of the numerator and the log of the denominator of
equation~(\ref{eq-bridged-lis1}) were found, and the overall standard
error was computed as the square root of the sum of the squares of
these two standard errors.  This method of converting estimates and
standard errors for $r$ to those for $\log(r)$ is valid
asymptotically.  It might be improved upon for finite samples, but
such improvements would probably not affect the relative merits of the
methods compared here.

\begin{figure}[p]

\centerline{\includegraphics{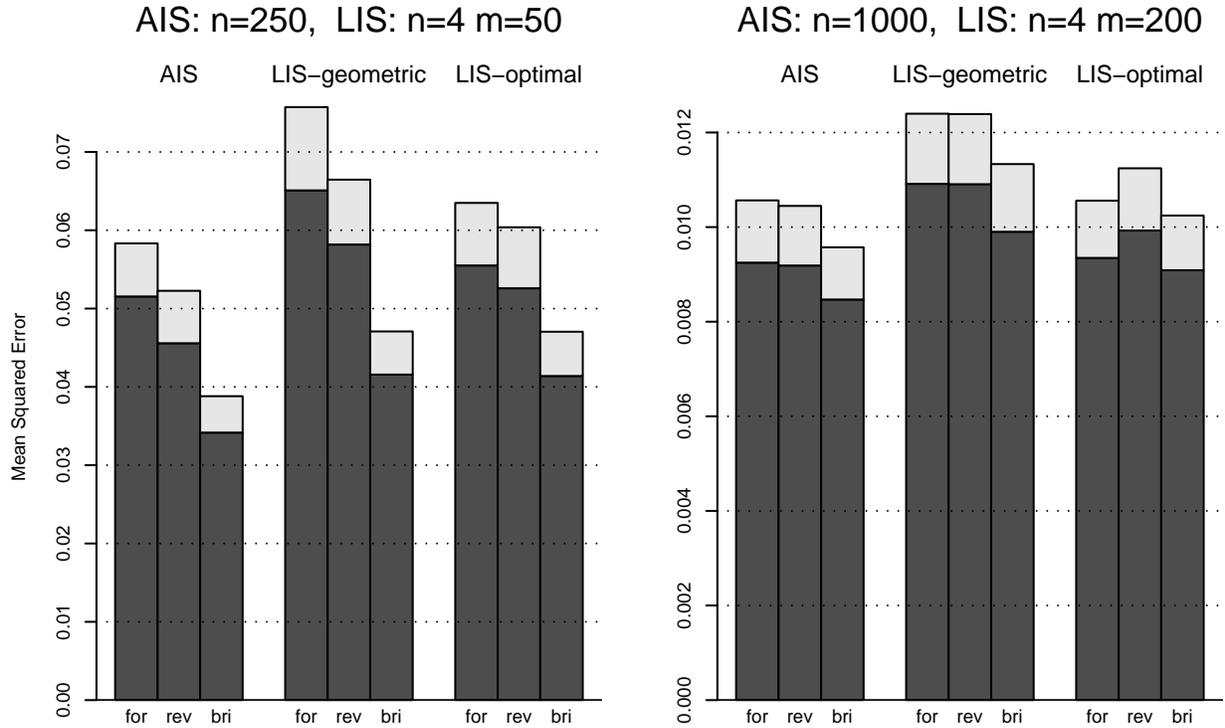}}

\vspace*{-8pt}


\caption[]{Results of short and long runs 
           on the distribution sequence with $s\!=\!1$, $t\!=\!4$, and 
           $q\!=\!2$.}\label{fig-r1}

\end{figure}

\begin{figure}[p]

\centerline{\includegraphics{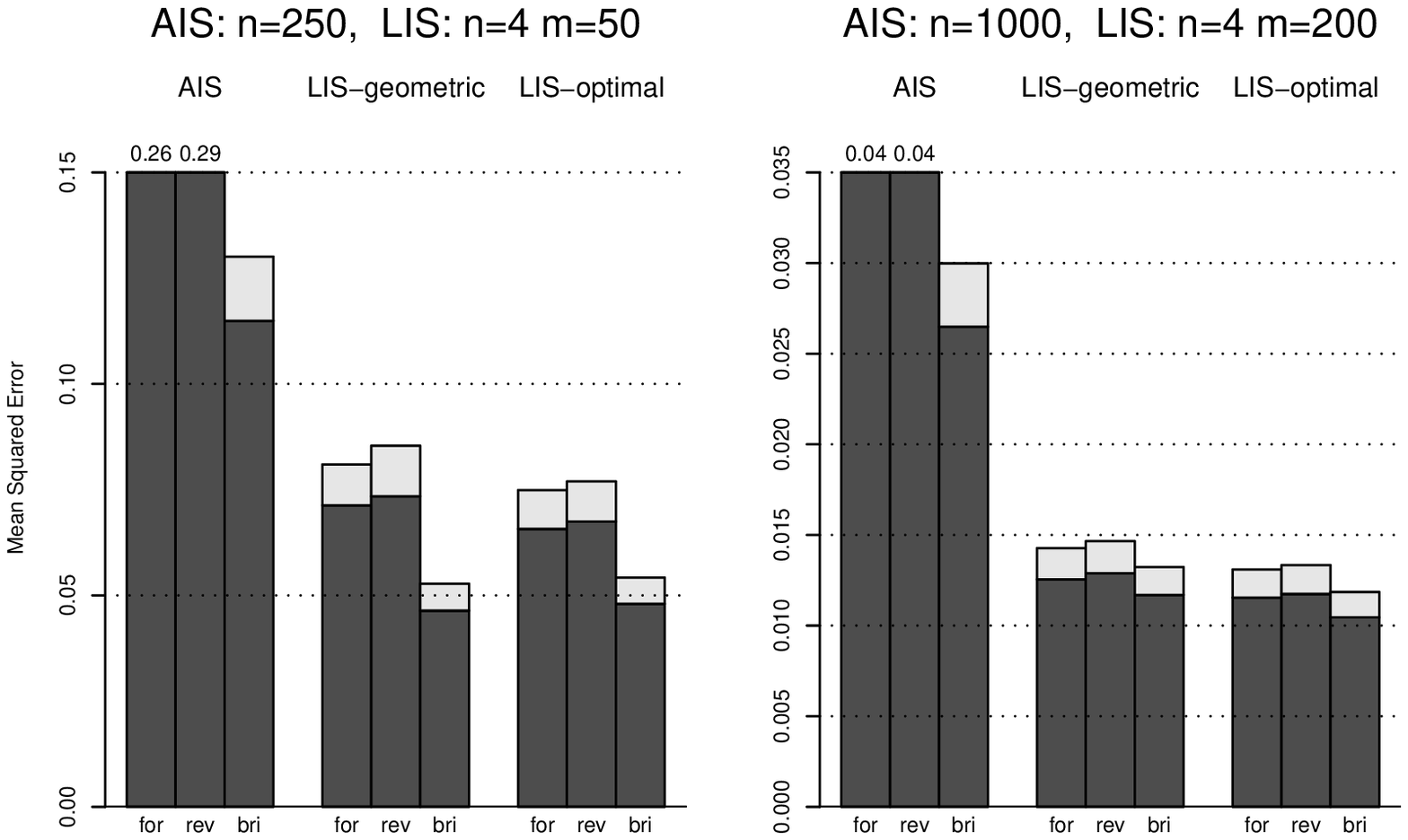}}

\vspace*{-8pt}


\caption[]{Results of short and long runs 
           on the distribution sequence with $s\!=\!1$, $t\!=\!4$, and 
           $q\!=\!10$.}\label{fig-r2}

\end{figure}

\begin{figure}[p]

\centerline{\includegraphics{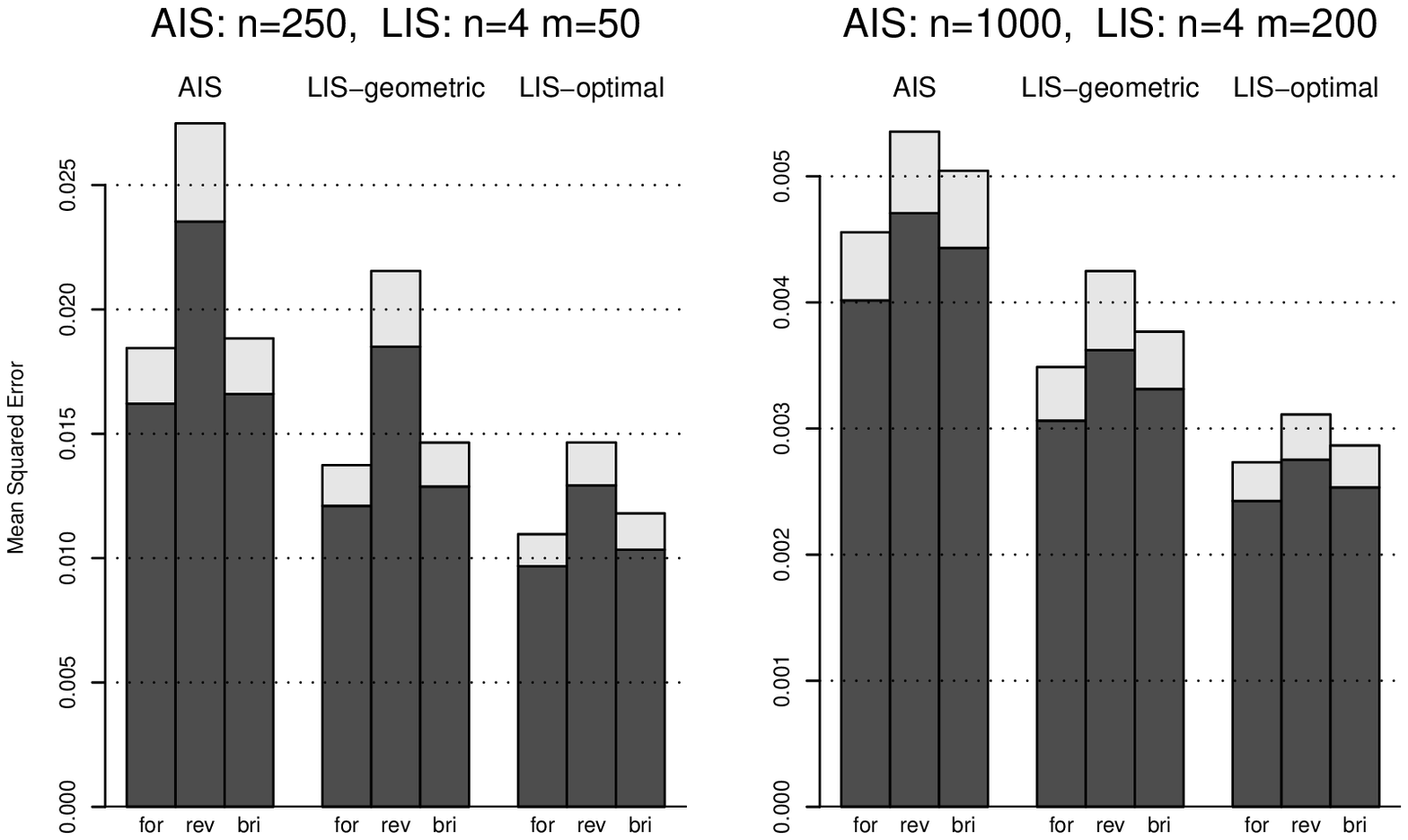}}

\vspace*{-8pt}


\caption[]{Results of short and long runs 
           on the distribution sequence with $s\!=\!0.05$, $t\!=\!0$, and 
           $q\!=\!2$.}\label{fig-r3}

\end{figure}

\begin{figure}[p]

\centerline{\includegraphics{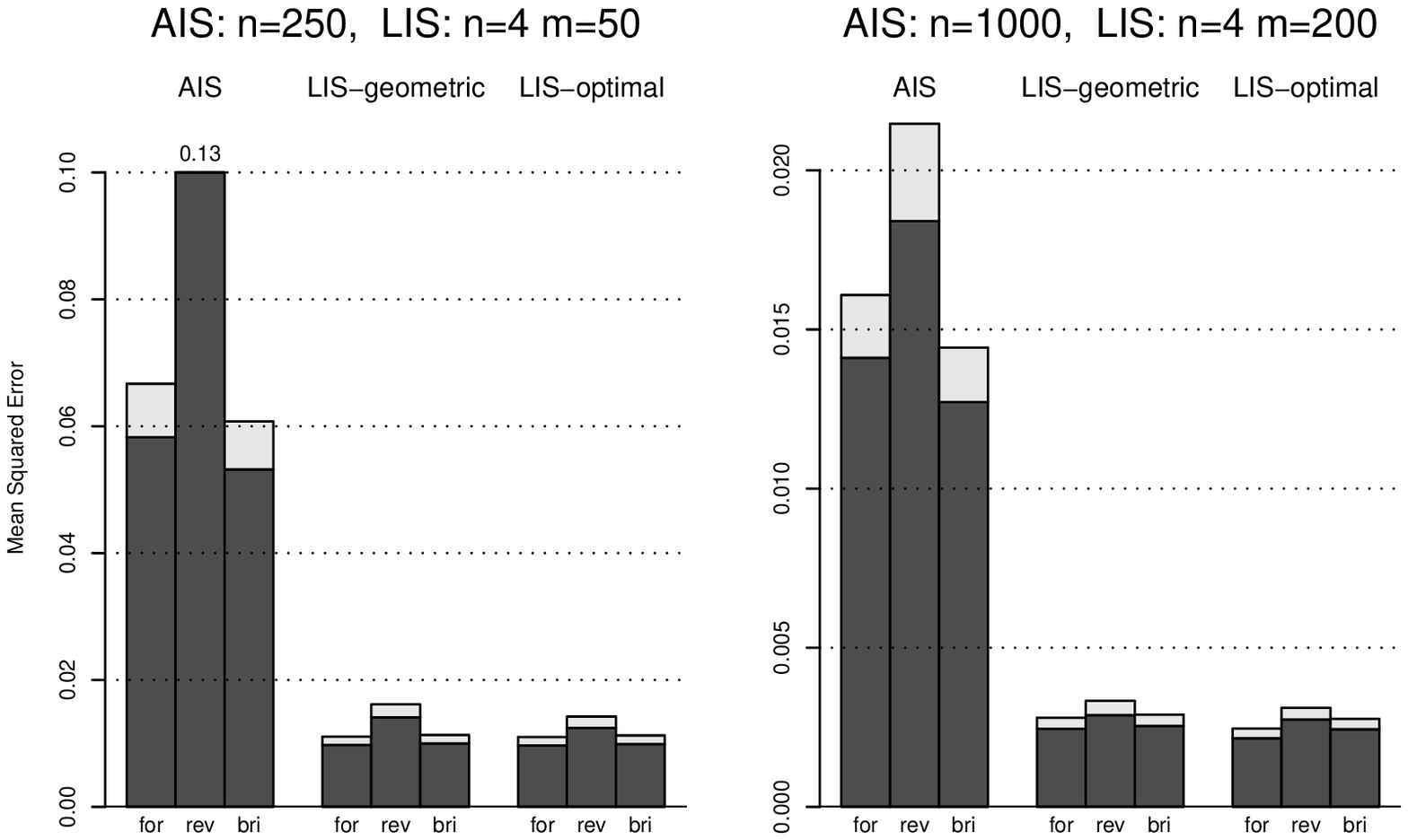}}

\vspace*{-8pt}

5\hspace*{0.1in}\makebox[3.2in]{Short Runs}\hfill\makebox[3.2in]{Long Runs}

\caption[]{Results of short and long runs 
           on the distribution sequence with $s\!=\!0.05$, $t\!=\!0$, and 
           $q\!=\!10$.}\label{fig-r4}

\end{figure}

\begin{figure}[p]

\centerline{\includegraphics{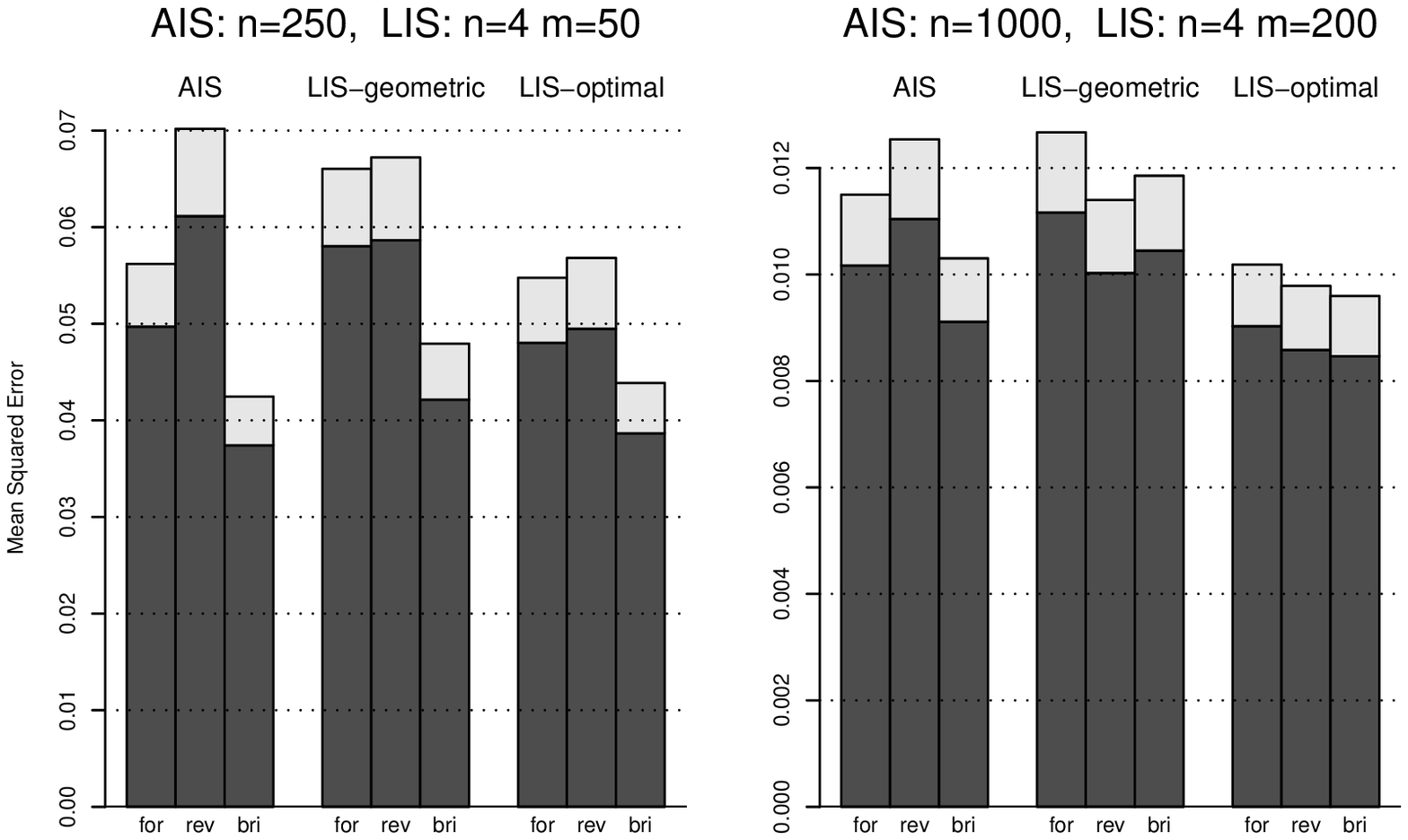}}

\vspace*{-8pt}


\caption[]{Results of short and long runs 
           on the distribution sequence with $s\!=\!0.3$, $t\!=\!2$, and 
           $q\!=\!2$.}\label{fig-r5}

\end{figure}

\begin{figure}[p]

\centerline{\includegraphics{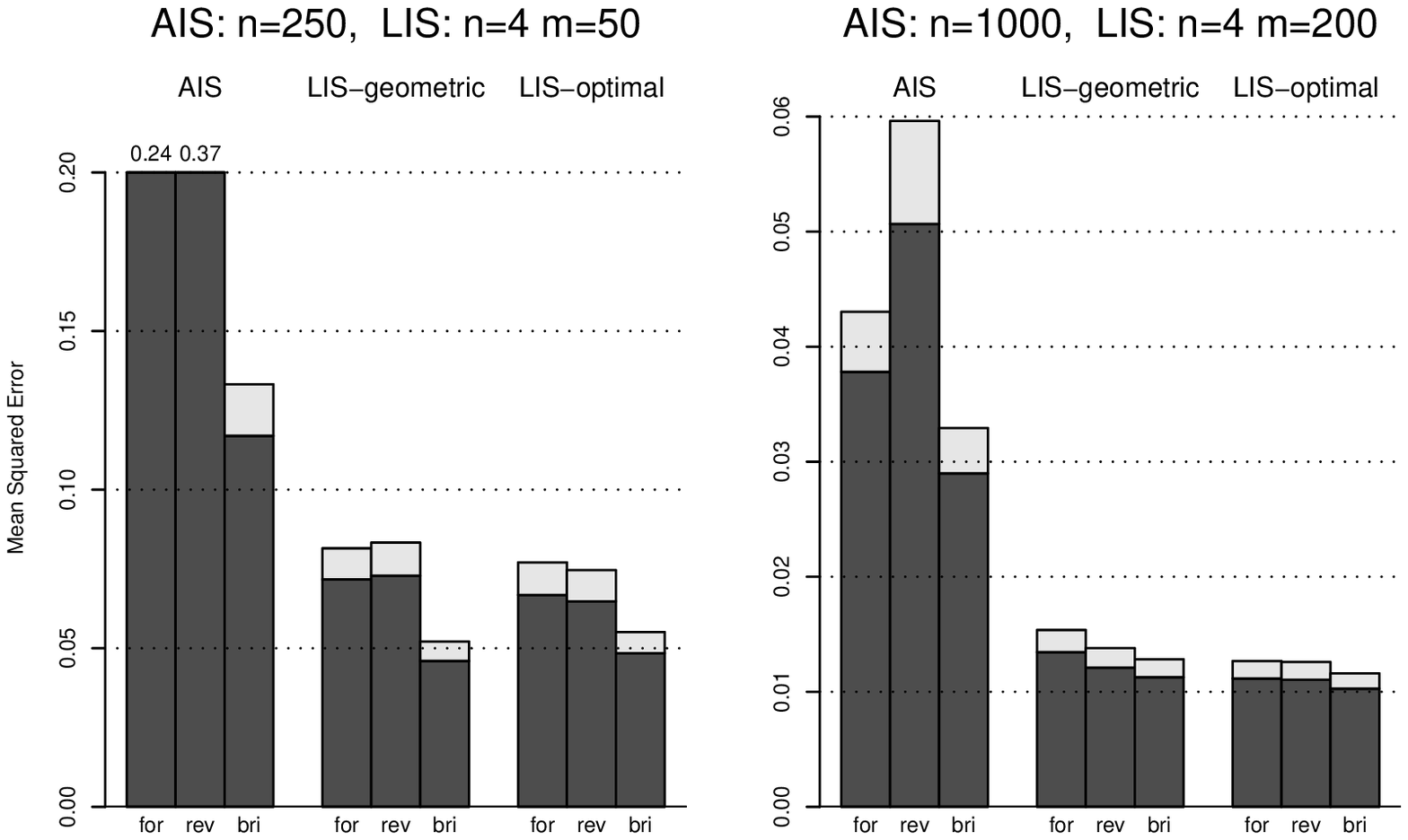}}

\vspace*{-8pt}


\caption[]{Results of short and long runs 
           on the distribution sequence with $s\!=\!0.3$, $t\!=\!2$, and 
           $q\!=\!10$.}\label{fig-r6}

\end{figure}

Figures~\ref{fig-r1} through \ref{fig-r6} plot the mean squared errors
of estimates for $\log(r)$ for the six sets of runs.  Results are
shown for AIS, for LIS using the geometric bridge, and for LIS using
the optimal bridge, with the true value of $r$.  Results for both the
forward and reverse versions of each method are shown, together with
the bridged version, using the optimal bridge, with $r$ obtained by
iteration.  Results for the short runs ($n\!=\!4$, $K_j\!=\!50$ for
LIS, $n\!=\!250$ for AIS) are on the left, and for the long runs
($n\!=\!4$, $K_j\!=\!200$ for LIS, $n\!=\!2000$ for AIS) on the right.
The mean squared error for each method was estimated by simulating
each method 2000 times, and comparing the estimates with the true
value of $\log(r)$.  The bars in the plots are dark up to the
estimated mean squared error minus twice its standard error, and are
then light up to the estimated mean squared error plus twice its
standard error.  For bars that extend above the plot the estimated
mean squared error is shown at the top of the bar.

The results for translated sequences of distributions ($t\!=\!4$ and
$s\!=\!1$) are shown in Figures~\ref{fig-r1} and~\ref{fig-r2}.  When the
distributions are Gaussian ($q\!=\!2$), no advantage is seen for LIS --- if
anything, LIS performs slightly worse than AIS, particularly when the
geometric bridge is used.  The forward and reverse forms of AIS and
LIS should have identical performance for these distribution
sequences, due to symmetry; any differences seen result from random
variation.  The bridged forms of both AIS and LIS perform better than
the unbridged forward and reverse forms.  The advantage of bridging is
less for the longer runs, however, as expected from the analysis at
the end of Section~\ref{sec-asym}.

When $q\!=\!10$, the distributions have much lighter tails than the
Gaussian, more closely resembling the uniform distributions analysed
in Section~\ref{sec-unif}.  For these sequences of distributions, LIS
performs substantially better than AIS.  The unbridged version of AIS
does particularly badly.  The mean squared error for the bridged
version of AIS is about 2.5 times greater than for the bridged version
of LIS.  It makes little difference whether the geometric or optimal
bridge is used for LIS.

Figures~\ref{fig-r3} and~\ref{fig-r4} show the results for sequences
of distributions with the same mean ($t\!=\!0$) but decreasing width
($s\!=\!0.05$).  For these sequences, a modest advantage of LIS over AIS
is apparent for the sequence of Gaussian distributions ($q\!=\!2$), with
the variance for AIS estimates being about a factor of 1.3 greater
than for LIS estimates with the geometric bridge, and about a factor
of 1.7 greater than for LIS estimates with the optimal bridge.  The
reversed AIS and LIS estimates are somewhat worse than the forward
estimates for this sequence of distributions.  No advantage is seen for
bridged AIS or LIS estimates.  

The results for the sequence of distributions with $q\!=\!10$ is similar,
except that the advantage of LIS over AIS is much greater --- about a
factor of 6.

Results for the last type of sequence, with $s\!=\!0.3$ and $t\!=\!2$, are
shown in Figures~\ref{fig-r5} and~\ref{fig-r6}.  This problem is a
hybrid of the previous two, with both translation and change in width,
producing results intermediate between those for the previous two
problems.  No difference in performance between AIS and LIS is
apparent for the Gaussian distributions ($q\!=\!2$), but the bridged forms
of both perform slightly better.  For the sequence of distributions
with $q\!=\!10$, a clear advantage of LIS over AIS can be seen, but this
advantage is not as great as for the sequence with $t\!=\!0$ and $s\!=\!0.05$.
The bridged forms of both AIS and LIS are again better, more so for
the short runs than for the long runs.

In addition to looking at the mean squared error of estimates found
with these methods, I also looked at the fraction of times that the
estimate for $\log(r)$ differed from the true value by more than twice
the standard error estimated using the $M$ runs.  This should be
approximately 5\% if the distribution of estimates is Gaussian, and
the standard errors are accurate.  For the longer runs, this fraction
was indeed near or only slightly above 5\% for all methods, except for
the unbridged AIS runs when these performed very poorly.  For the
shorter runs, however, the unbridged AIS and LIS methods produced
estimates more than two standard errors from the mean around 10\% of
the time (sometimes much more often, when unbridged AIS performed
poorly).  Both the bridged AIS and the bridged LIS methods gave more
reliable standard errors.  However, it is possible that better
standard errors for the unbridged methods might be obtained with a
more sophisticated approach than I used.

\begin{figure}[p]

\centerline{\includegraphics{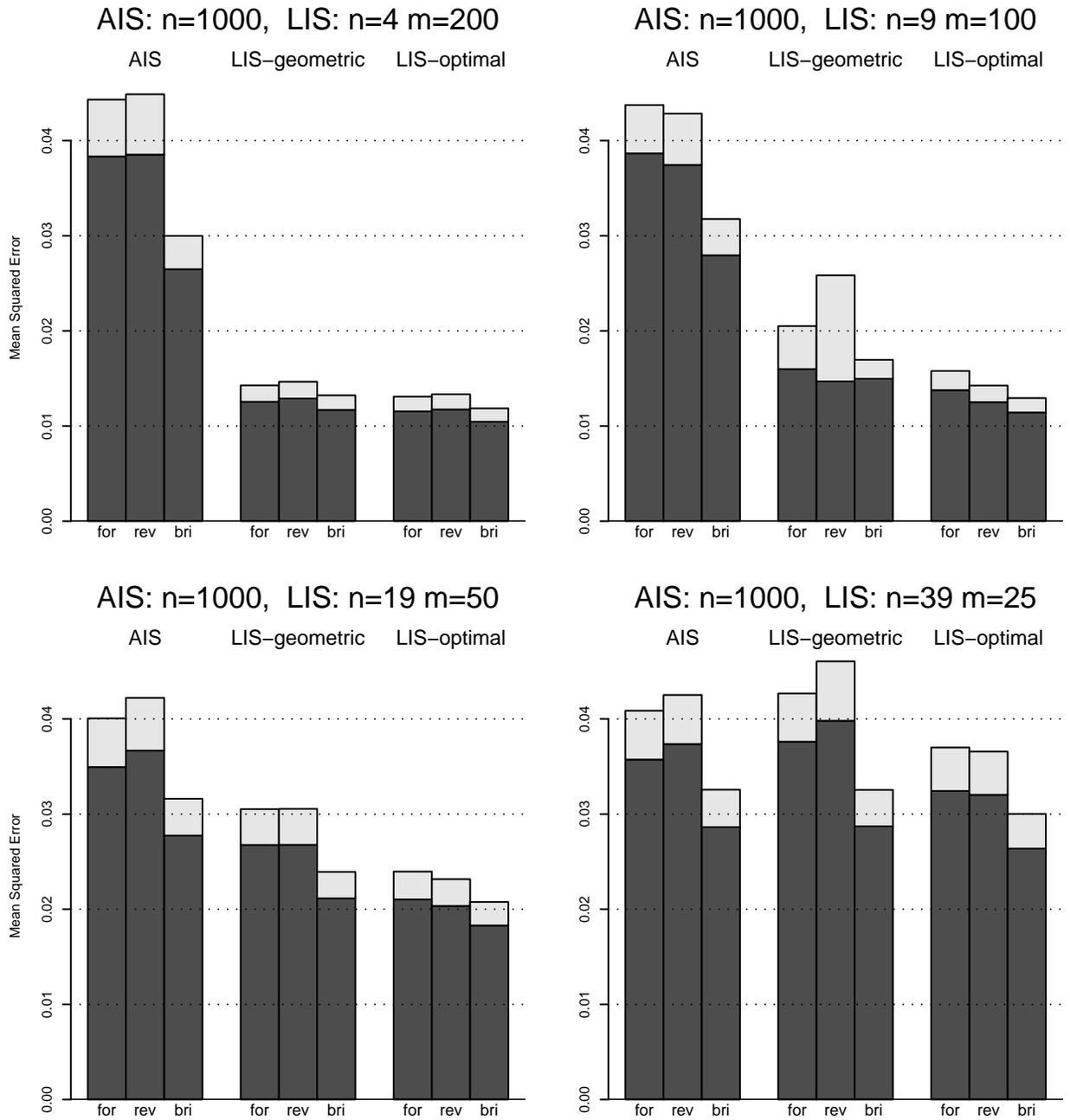}}

\vspace*{-8pt}

\caption[]{Results using increasing values of $n$ for LIS, while keeping
           computation time constant, for the distribution sequence with
           $s\!=\!1$, $t\!=\!4$, and $q\!=\!10$.  The same AIS procedure was 
           used for all plots, but results vary randomly.}\label{fig-r7}

\end{figure}

\begin{figure}[p]

\centerline{\includegraphics{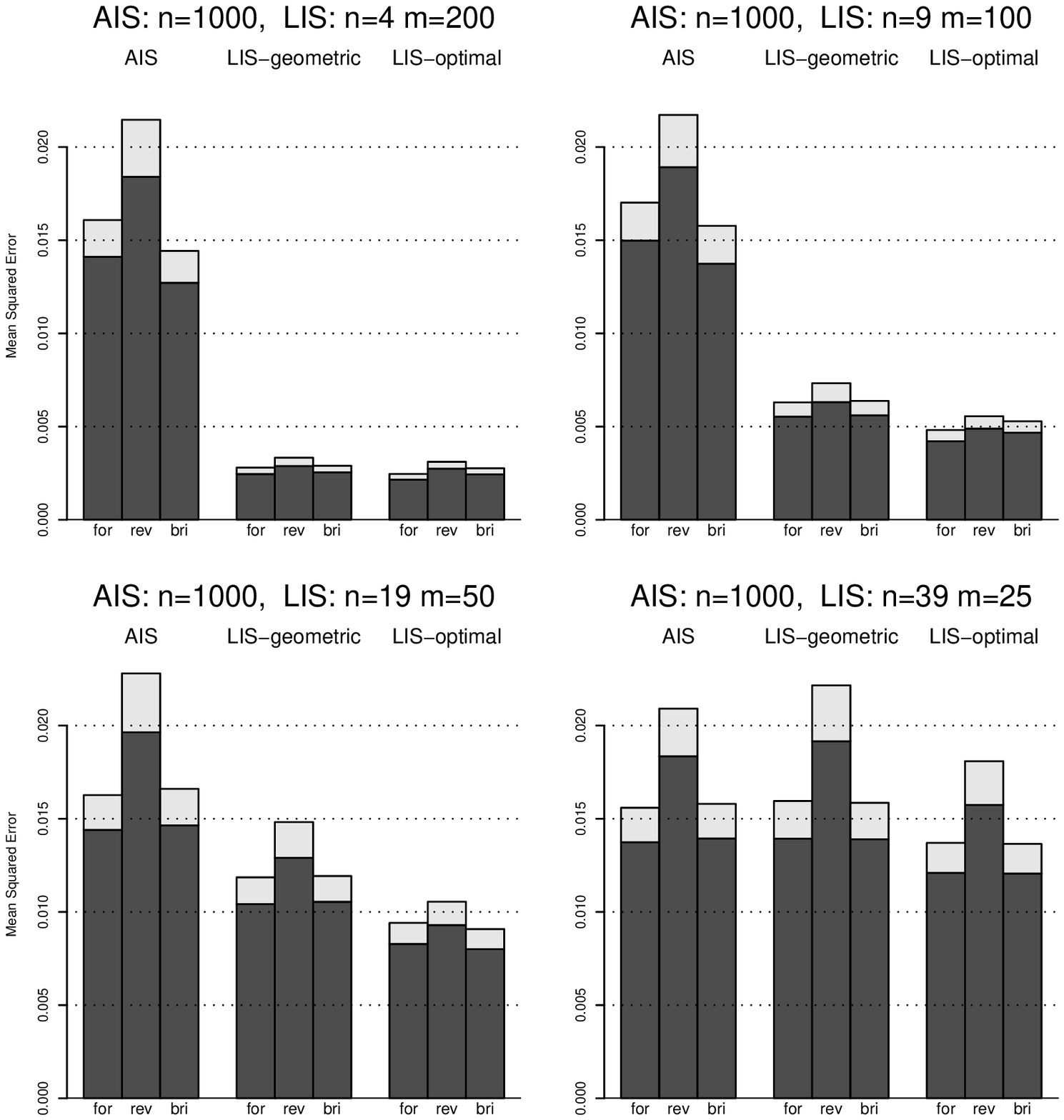}}

\vspace*{-8pt}

\caption[]{Results using increasing values of $n$ for LIS, while keeping
           computation time constant, for the distribution sequence with
           $s\!=\!0.05$, $t\!=\!0$, and $q\!=\!10$.  The same AIS procedure was
           used for all plots, but results vary randomly.}\label{fig-r8}

\end{figure}

I performed additional runs to verify and extend some of the analytic
results from Section~\ref{sec-anal}.  Figures~\ref{fig-r7}
and~\ref{fig-r8} show results obtained using LIS with increasing
numbers of intermediate distributions, starting with the value of
$n\!=\!4$ used for the tests above, and continuing to $n\!=\!9$,
$n\!=\!19$, and $n\!=\!39$, while keeping the computation time
constant by decreasing $m$ in proportion to $n\!+\!1$.  The two
distribution sequences with $s\!=\!1$ and $t\!=\!4$ and with
$s\!=\!0.05$ and $t\!=\!0$ were used, in both cases with $q\!=\!10$.
The sequence with $t\!=\!0$ and $s\!=\!0.05$ has the form of
equation~(\ref{eq-U-dist}), so in accordance with the analysis of
Section~\ref{sec-asym}, we expect that asymptotically, as $n$
increases, LIS and AIS should have the same performance.  This is
indeed what we see in Figure~\ref{fig-r7}.  We also see the same
behaviour for the sequence with $t\!=\!4$ and $s\!=\!1$ in
Figure~\ref{fig-r8}.

\begin{figure}[p]

\centerline{\includegraphics{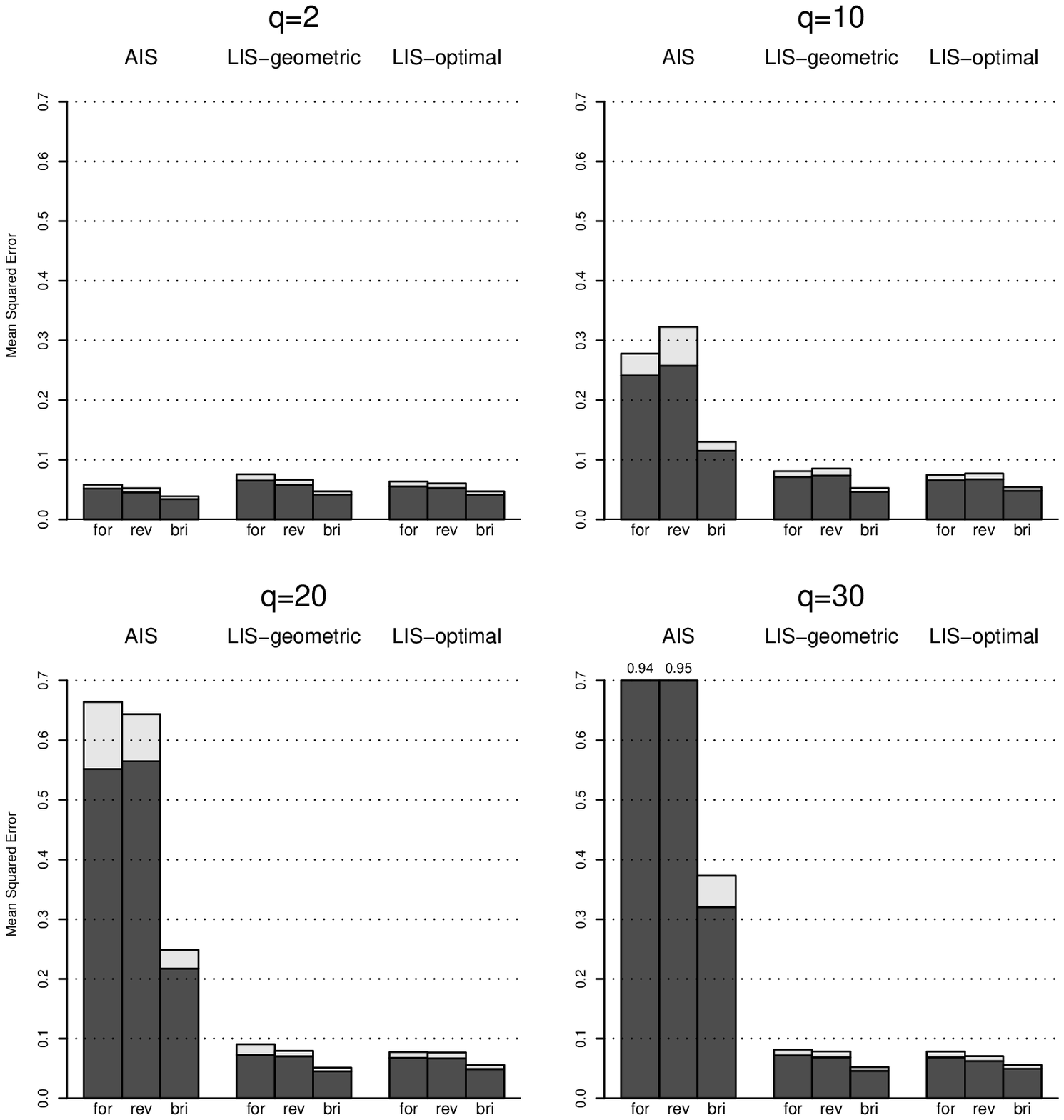}}

\vspace*{-8pt}

\caption[]{Results with increasing values of $q$, for sequences of
           distributions with $s\!=\!1$ and $t\!=\!4$.  The AIS runs used
           $n\!=\!250$; the LIS runs used $n\!=\!4$ and $m\!=\!50$,
           requiring the same amount of computation.}\label{fig-r9}

\end{figure}

\begin{figure}[p]

\centerline{\includegraphics{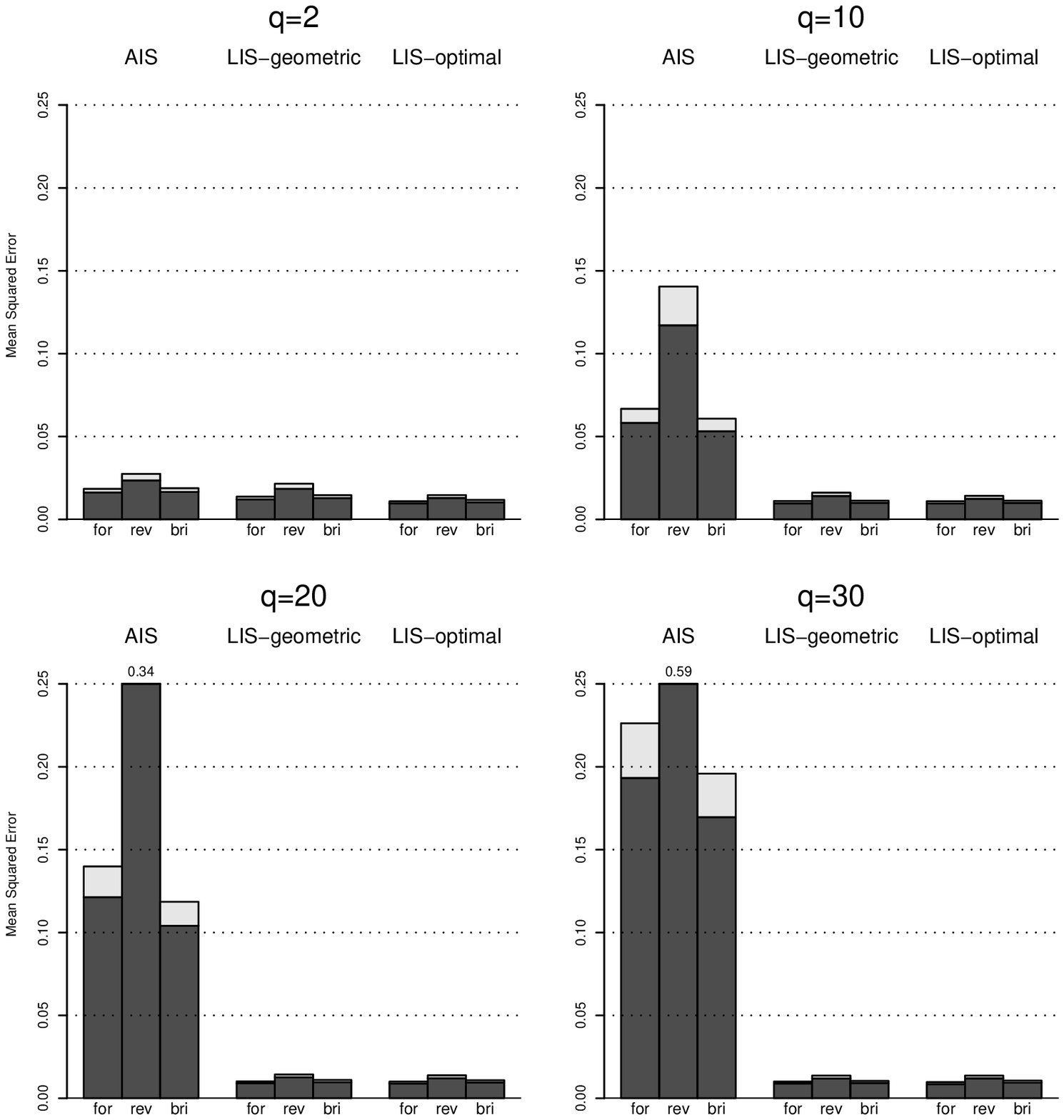}}

\vspace*{-8pt}

\caption[]{Results with increasing values of $q$, for sequences of
           distributions with $s\!=\!0.05$ and $t\!=\!1$.  The AIS runs used
           $n\!=\!250$; the LIS runs used $n\!=\!4$ and $m\!=\!50$, 
           requiring the same amount of computation.}\label{fig-r10}

\end{figure}

As $q$ increases, the distributions become close to uniform, and the
results of Section~\ref{sec-unif} should apply.  To test this, I tried
values of $q\!=\!2$, $q\!=\!10$, $q\!=\!20$, and $q\!=\!30$ for the
distribution sequence with $s\!=\!1$ and $t\!=\!4$ and the sequence with
$s\!=\!0.05$ and $t\!=\!0$.  Results are shown in Figures~\ref{fig-r9}
and~\ref{fig-r10}.  (The results for $q\!=\!2$ and $q\!=\!10$ are the same as
on the left in Figures~\ref{fig-r1} to~\ref{fig-r4}, though the scale
differs.)  

For the sequences with $s\!=\!1$ and $t\!=\!4$, the limiting uniform
distributions have the form of the second example in
Section~\ref{sec-unif}.  As noted there, AIS estimates do not converge
to the correct value of $r$ for this distribution sequence; bridged AIS
estimates do converge, but may be rather inefficient.  We see
analogous behaviour in Figure~\ref{fig-r9} when $q$ is large.  The
mean squared error of the AIS estimates increases approximately
linearly with $q$ over the range $q\!=\!10$ to $q\!=\!30$.  The
bridged AIS estimates also get worse as $q$ increases, but more
slowly.  In contrast, the mean squared error of the LIS estimates
changes hardly at all as $q$ increases.

The story is similar for sequences with $s\!=\!0.05$ and $t\!=\!1$,
for which the limiting uniform distributions correspond to those in
the first example of Section~\ref{sec-unif}.  The LIS estimates
perform about equally well for all values of $q$, but the AIS
estimates are dramatically worse for large values of $q$.  For this
sequence, reverse AIS estimates are much worse than forward AIS
estimates, and bridging does not help.

According to the analysis of Section~\ref{sec-asym}, the choice of
choice of $n\!=\!4$ for LIS used above is not optimal for either of
these distribution sequences when $q$ is large.  For the sequence with
$s\!=\!1$ and $t\!=\!4$, using $n\!=\!6$ should be better by a factor
of 1.176.  However, in LIS runs with $q=30$, the mean squared error
using $n=\!=\!4$ and $m\!=\!200$ is indistinguishable from that using
$n\!=\!6$ and $m\!=\!143$, given the standard errors (a factor of 1.09
or more should have been detectable).  Of course, $q=30$ does not give
exactly uniform distributions, and these values of $m$ may not be
large enough for the asymptotic results to apply, especially since the
Markov transitions do not sample independently.  For the sequence with
$s\!=\!0.05$ and $t\!=\!0$, the results in Section~\ref{sec-asym}
indicate that using $n\!=\!3$ should be better by a factor of 1.084.
In this case, LIS runs with $q=30$ using $n\!=\!3$ and $m\!=\!250$ are
better than runs using $n\!=\!4$ and $m\!=\!200$ by a factor of 1.16,
significantly greater than one given the standard errors, but not
significantly different from the expected ratio of 1.084.

\section{\hspace*{-7pt}Other applications of linked 
                       sampling}\label{sec-gen}\vspace*{-10pt}

So far in this paper, I have focused on how Linked Importance Sampling
can be used to estimate ratios of normalizing constants.  LIS can also
be used to estimate expectations with respect to $\pi_1$, however, and
in some applications, this may be its most important use.  Linked
sampling methods related to LIS can also be applied in other ways.  I
briefly described these other applications here, outlining the use of
linked sampling for `dragging' fast variables in some detail.

\subsection{\hspace*{-4pt}Estimating expectations}\vspace*{-4pt}

The expectation of some function, $a(x)$, with respect to $\pi_1$ 
can be estimated using simple importance sampling, with points drawn
from $\pi_0$, as follows:
\beq
   E_{\pi_1}\big[a(X)\big] 
   \ \ = \ \ E_{\pi_0}\!\left[ a(X) {p_1(X) \over p_0(X)}\right] \, \Big/\ 
             {Z_1 \over Z_0}
   \ \ \approx\ \ 
{1 \over N}\sum_{i=1}^N\, a(x^{(i)})\, {p_1(x^{(i)}) \over p_0(x^{(i)})}\ \Big/\
{1 \over N}\sum_{i=1}^N\, {p_1(x^{(i)}) \over p_0(x^{(i)})}
\label{eq-is-exp}
\eeq
where $x^{(i)},\ldots,x^{(N)}$ are drawn from $\pi_0$.
Like equation~(\ref{eq-simple}), this estimate is valid only if 
no region having zero probability under $\pi_0$ has non-zero probability 
under $\pi_1$.  The two factors of $1/N$ of course cancel, but are included 
to emphasize the connection with the estimate for $r=Z_1/Z_0$, which is
simply the denominator of the estimate above.

Since LIS can be viewed as simple importance sampling on an extended
state space, with distributions $\Pi_0$ and $\Pi_1$ defined by the
forward and reverse procedures of Section~\ref{sec-lis}, we can use
equation~(\ref{eq-is-exp}) to estimate any quantity that can be
expressed as an expectation with respect ot $\Pi_1$.  Step (1) of the
reverse procedure defining $\Pi_1$ sets $x_{n,\mu_n}$ to a value
randomly chosen from $\pi_{\eta_n} = \pi_1$.  Step (2) then sets the
other $x_{n,k}$ to values obtained from $x_{n,\mu_n}$ by applying
Markov chain transitions that leave $\pi_1$ invariant.  It follows
that under $\Pi_1$, all the points $x_{n,k}$ have marginal
distribution $\pi_1$ (though they may not be independent).  Accordingly,
\beq
  E_{\pi_1}\big[a(X)\big] & = & E_{\,\Pi_1}\!\left[ 
    {1 \over K_n\!+\!1}\, \sum_{k=0}^{K_n} a(X_{n,k}) \right]
\eeq
Estimating the right side as in equation~(\ref{eq-is-exp}), and using
the fact that the ratio of probabilities under $\Pi_1$ over those
under $\Pi_0$ is given by $\rhatLIS^{(i)}$ in equation~(\ref{eq-lis}),
we get the estimate
\beq
  E_{\pi_1}\big[a(X)\big] & \approx & 
  \sum_{i=1}^M {\rhatLIS^{(i)} \over K_n\!+\!1} 
     \sum_{k=0}^{K_n} a(x^{(i)}_{n,k}) 
  \ \Big/\ 
  \sum_{i=1}^M \rhatLIS^{(i)}
\label{eq-lis-exp}
\eeq

If the $M$ runs of LIS are started by sampling independently from
$\pi_0$ (as will often be possible), the standard error of this
estimate can be assessed in the usual fashion for importance sampling,
as I have discussed for the analogous AIS estimates in (Neal 2001).
This error assessment can be difficult, since when some
$\rhatLIS^{(i)}$ are much larger than others, the variance of
$\rhatLIS^{(i)}$ is hard to estimate.  Note, however, that the degree
to which the Markov chain transitions used have converged need not be
assessed, a possible advantage compared with simple MCMC estimates.  The
estimate of equation~(\ref{eq-lis-exp}) will be asymptotically correct
(as $M\rightarrow\infty$) regardless of how far these Markov chain
transitions are from convergence.

The primary reason one might wish to use LIS to estimate expectations
is that going through the sequence of distributions parameterized by
$\eta_0,\ldots,\eta_n$ may produce an `annealing' effect, which
prevents the Markov chain sampler from being trapped in a local mode
of the distribution.  Compared with the analogous AIS procedure, LIS
may perform better for some forms of distributions, for the same
reasons as were discussed in Sections~\ref{sec-anal}
and~\ref{sec-cmp}.  One should also note that LIS estimates for
expectations with respect to $\pi_{\eta_j}$ for all $j$ can easily be
obtained from a single set of runs, by simply considering the results
of each LIS run up to the point where the sample for $\pi_{\eta_j}$ is
obtained.

\subsection{\hspace*{-4pt}A linked form of tempered transitions}\vspace*{-4pt}

My `tempered transition' method (Neal 1996) is another approach to
sampling from distributions with isolated modes, between which
movement is difficult for Markov chain transitions such as simple
Metropolis updates.  In this approach, such simple Markov chain
transitions are supplemented by occasional complex `tempered
transitions', composed of many simple Markov chain transitions.  A
tempered transition consists of several stages, which proceed through
a sequence of distributions, from the distribution being sampled, to a
`higher temperature' distribution in which movement between modes is
easier, and then back down to the distribution being sampled.  At each
stage of a tempered transition, we generate a single new state by
applying a Markov chain transition to the current state, after which
we switch to the next distribution in the sequence. The second half of
a tempered transition is similar to an Annealed Importance Sampling
run, while the first half is similar to an AIS run with the reversed
sequence of distributions.

A similar `linked' procedure can be defined, in which at each stage we
generate a chain of states by applying a Markov chain transition.
We then select a `link state' from this sequence (using a suitable
bridge distribution) which serves as the starting point for the chain
of states generated in the next stage.  In the final stage, a chain of
states is produced using a Markov chain transition that leaves the
distribution being sampled invariant, and a candidate state is
selected uniformly at random from this chain.  The appropriate
probability for accepting this candidate state is computed using
ratios similar to those going into the LIS estimate of
equation~(\ref{eq-lis}).

As discussed in Section~\ref{sec-cmp}, for AIS to work well, all
distributions in the sequence must assign reasonably high probability
to regions of the space that have non-negligible probability under the
next distribution in the sequence.  One would expect tempered
transitions to work well only when this holds for both the sequence
and its reversal.  In contrast, one would expect the `linked' version
of tempered transitions to work well as long as the sequence satisfies
the weaker condition that there be some `overlap' between adjacent
distributions (assuming a suitable bridge distribution is used).

\subsection{\hspace*{-4pt}Dragging fast variables using linked 
                          chains}\vspace*{-4pt}

A slight modification of the tempered transition method can be applied
to problems in which the state is composed of both `fast' and `slow'
variables.  We will write the distribution of interest for such a problem
as 
\beq
   \pi(x,y) & = & (1/Z)\, \exp(-U(x,y)) 
\eeq 
where $x$ denotes the `fast' variables and $y$ the `slow' variables.
We assume that the computation is dominated by the time required to
evaluate $U(x,y)$, but that once $U(x,y)$ has been evaluated, with
relevant intermediate quantities saved,
evaluating $U(x',y)$ for any new $x'$ is much faster than evaluating
$U(x',y')$ for some $y'$ not previously encountered.  One example of
such a problem is inference for Gaussian process classification models
(Neal 1999), in which $y$ consists of the hyperparameters defining the
covariance function used, and $x$ consists of the latent variables
associated with the $n$ observations.  After a change to $y$, we must
recompute the Cholesky decomposition of an $n \times n$ covariance matrix, 
which takes time proportional to $n^3$, whereas after a change to $x$
only, $U(x,y)$ can be re-computed in time proportional to $n^2$,
assuming the Cholesky decomposition for this value of $y$ has been
saved.

In my method for `dragging' fast variables (Neal 2004), the ability
to quickly re-evaluate $U(x,y)$ when only $x$ changes is exploited to
allow larger changes to be made to $y$ than would be possible if $x$
were kept fixed, or were given a new value from some simple proposal
distribution.  From the state $(x_0,y_0)$, a dragging
update proposes a new value $y_1$, drawn from some symmetrical proposal
distribution, in conjunction with a new value $x_1$ that is found by 
applying a succession of Markov chain updates that leave
invariant distributions in the series, $\pi_{\eta_j}(x)$, for 
$j=1,\ldots,n\!-\!1$, with $0<\eta_j<\eta_{j+1}<1$.  The proposed state,
$(x_1,y_1)$, is then accepted or rejected in a fashion analogous to tempered 
transitions.  

The distributions in the sequence used are defined by the following
unnormalized probability or density function, which depends on the
current and proposed values for $y$:
\beq
  p_{\eta}(x) & = & 
    \exp\,(\,-\,((1\!-\!\eta)\, U(x,y_0)\ +\ \eta\, U(x,y_1)))
\label{eq-drag-p}
\eeq
The corresponding normalized probability or density function will be
written as $\pi_{\eta}$.  Note that $\pi_0(x) = \pi(x|y_0)$ and 
$\pi_1(x)=\pi(x|y_1)$.  Crucially,
after $U(x,y_0)$ and $U(x,y_1)$ have been evaluated once (for any~$x$),
we can evaluate $p_{\eta}(x)$ for any $\eta$ and any $x$
without any further `slow' computations.
Indeed, since $U(x_0,y_0)$ will usually have already been evaluated as part of 
the previous Markov chain transition, only one slow computation will be required
to evaluate $p_{\eta}(x)$ for any number of values of $\eta$ and $x$.

A `linked' dragging update can be defined as follows.  Given
the sequence of distributions defined by $\eta_0,\ldots,\eta_n$, with 
$\eta_0=0$ and $\eta_n=1$, the numbers of transitions ($T$ or $\underline{T}$)
to perform for each distribution over $x$, denoted by $K_0,\ldots,K_n$, and a 
set of bridge distributions, denoted by $p_{j*j+1}$, for $j=0,\ldots,n\!-\!1$, 
an update from the current state $(x_0,y_0)$ is done as follows:\vspace*{5pt}

\begin{center}\bf The Linked Dragging Procedure\end{center}\vspace*{-5pt}

\begin{enumerate}
\item[1)] Propose a new value, $y_1$, from some proposal distribution
          $S(y_1|y_0)$, which satisfies the symmetry condition that $S(y_1|y_0)
          =S(y_0|y_1)$.
\item[2)] Pick an integer $\nu_0$ uniformly at random from $\{0,\ldots,K_0\}$,
          and then set $x_{0,\nu_0}$ to the current values of the fast
          variables, $x_0$.
\item[3)] For $j\,=\,0,\ldots,n$, create a chain of values for $x$ associated 
          with $\pi_{\eta_j}$ as follows:
\begin{enumerate}
  \item[a)] If $j>0$:\ \ Pick an integer $\nu_j$ uniformly at random from 
            $\{0,\ldots,K_j\}$, and then set $x_{j,\nu_j}$ to $x_{j-1*j}$.
  \item[b)] For $k\,=\,\nu_j+1,\ldots,K_j$, draw $x_{j,k}$ according to the
            forward Markov chain transition probabilities 
            $T_{\eta_j}(x_{j,k-1},x_{j,k})$.  (If $\nu_j=K_j$, do nothing in 
            this step.)
  \item[c)] For $k\,=\,\nu_j-1,\ldots,0$, draw $x_{j,k}$ according to the 
            reverse Markov chain transition probabilities
            $\underline{T}_{\eta_j}(x_{j,k+1},x_{j,k})$. (If $\nu_j=0$, do 
            nothing in this step.)
  \item[d)] If $j<n$:\ \ Pick a value for $\mu_j$ from
            $\{0,\ldots,K_j\}$ according to the following probabilities
            \beq
              \Pi_0(\mu_j\,|\,x_j) & = & 
                 {p_{j*j+1}(x_{j,\mu_j}) \over p_{\eta_j}(x_{j,\mu_j})}
                 \ \Big/\
                 \sum_{k=0}^{K_j} {p_{j*j+1}(x_{j,k}) \over p_{\eta_j}(x_{j,k})}
            \eeq
            and then set $x_{j*j+1}$ to $x_{j,\mu_j}$.
\end{enumerate}
\item[3)] Set $\mu_n$ to a value chosen uniformly at random from 
          $\{0,\ldots,K_n\}$, and let the proposed new values for the fast
          variables, $x_1$, be equal to $x_{n,\mu_n}$.
\item[4)] Accept $(x_1,y_1)$ as the new state with probability
\beq
   \min \left\{\, 1,\ \
   \prod_{j=0}^{n-1} \left[
     {1 \over K_j+1}\, \sum_{k=0}^{K_j}\,
          { p_{j*j+1}(x_{j,k}) \over p_{\eta_j}(x_{j,k}) }
     \ \Big/\
     {1 \over K_{j+1}+1}\, \sum_{k=0}^{K_{j+1}}\, 
          { p_{j*j+1}(x_{j+1,k}) \over p_{\eta_{j+1}}(x_{j+1,k}) }
     \right]
   \,\right\}
\eeq
          If $(x_1,y_1)$ is not accepted, the new state is the same as
          the old state, $(x_0,y_0)$.\vspace*{-6pt}
\end{enumerate}
One can show that this update leaves $\pi(x,y)$ invariant by showing
that it satisfies detailed balance, which in turns follows from the
stronger property that the probability of starting at $(x_0,y_0)$,
assuming this start state comes from $\pi(x,y)$, then generating the various
quantities produced by the above procedure, and finally accepting $(x_1,y_1)$
as the new state, is the same as the probability of starting this procedure
at $(x_1,y_1)$, generating the same quantities in reverse, and finally accepting
$(x_0,y_0)$. The proof of this is analogous to the derivation of LIS in 
Section~\ref{sec-lis}.

To use the linked dragging procedure, we need to select suitable
bridge distributions.  Since the characteristics of $\pi_{\eta}(x)$
will depend on $y_0$ and $y_1$, and of course $\eta$, we may not know
enough to select good estimates for the values of $r$ needed to use
the optimal bridge of equation~(\ref{eq-opt-bridge}), though we might
try just setting $r$ to one.  This is not a problem for the geometric bridge of 
equation~(\ref{eq-geo-bridge}), for which the acceptance probability
above can be written as\vspace*{2pt}
\beq
   \min \left\{\, 1,\ \
   \prod_{j=0}^{n-1} \left[
     {1 \over K_j+1}\, \sum_{k=0}^{K_j}\,
          \sqrt{{ p_{\eta_{j+1}}(x_{j,k}) \over p_{\eta_j}(x_{j,k}) }}
     \ \Big/\
     {1 \over K_{j+1}+1}\, \sum_{k=0}^{K_{j+1}}\, 
          \sqrt{{ p_{\eta_j}(x_{j+1,k}) \over p_{\eta_{j+1}}(x_{j+1,k}) }}\,
     \right]
   \,\right\}\\[-10pt]\nonumber
\eeq
From equation~(\ref{eq-drag-p}), we see that 
\beq
  { p_{\eta_{j+1}}(x_{j,k}) \over p_{\eta_j}(x_{j,k}) }
  & = & \exp\,(\,-\,(\eta_{j+1}\!-\!\eta_j)\,
                    (U(x_{j,k},y_1)\!-\!U(x_{j,k},y_0))) \\[6pt]
  { p_{\eta_j}(x_{j+1,k}) \over p_{\eta_{j+1}}(x_{j+1,k}) }
  & = & \exp\,(\,-\,(\eta_{j+1}\!-\!\eta_j)\,
                    (U(x_{j+1,k},y_0)\!-\!U(x_{j+1,k},y_1))) 
\eeq
For the simplest case with no intermediate distributions (ie, with $n\!=\!1$), 
the acceptance probability simplifies to
\beq
   \min \left\{\, 1,\ \
   { \displaystyle {1 \over K_0+1}\, \sum_{k=0}^{K_0}\,
       \exp\,(\,-\,(U(x_{j,k},y_1)\!-\!U(x_{j,k},y_0))\,/\,2)
     \over
     \displaystyle {1 \over K_1+1}\, \sum_{k=0}^{K_1}\, 
       \exp\,(\,-\,(U(x_{j,k},y_0)\!-\!U(x_{j,k},y_1))\,/\,2)
   } \right\}
\eeq

\section{\hspace*{-7pt}Conclusions and Future work}\vspace*{-10pt}

In this paper, I have demonstrated that in some situations Linked
Importance Sampling is substantially more efficient than Annealed
Importance Sampling, provided a suitable number of intermediate
distributions are used.  However, in other situations, where the tails
of the distributions involved are sufficiently heavy, the two methods
are about equally efficient.  More research is therefore needed to
determine for which problems of practical interest LIS, and related
linked sampling methods, will be useful.

In tests on multivariate Gaussian distributions, I have not seen an
advantage for LIS over AIS.  Both perform about equally well on a
sequence of 100-dimensional spherical Gaussian distributions with
variances changing by a factor of two, so that $\log(r) = -100$.  This
is in accord with the results in Section~\ref{sec-cmp}, where LIS had
little or no advantage over AIS when the distributions were Gaussian.
LIS is more likely to be useful for problems involving continuous
distributions with lighter tails.

One problem that may benefit from LIS is that of computing the
probability of a very rare event, which can be cast as computing the
normalizing constant for a distribution with the constraint that the
state be in the set corresponding to this event.  Intermediate
distributions might use looser forms of this constraint.  If, in all
these distributions, states violating the constraints have zero
probability, AIS will tend to have the same bad behaviour seen with
uniform distributions in Section~\ref{sec-unif}, while LIS may work
much better.

Another context where LIS may outperform AIS is when only a fixed
number of intermediate distributions are available --- ie, only a
finite number of values are allowed for $\eta$.  This is the situation
for the `sequential importance sampler' of MacEachern, Clyde, and Liu
(1999), which can be seen as an instance of AIS (Neal 2001).  Here,
the intermediate distributions use only a fraction of the $n$ items in
the data set; such a fraction can only have the form $j/n$ with $j$ an
integer.  The distance between successive distributions for this
problem may sometimes be too great for AIS to work well, but their
overlap might nevertheless be sufficient for LIS.

It may be possible to improve LIS by reducing the variance in how well
it samples at each stage.  Instead of performing a predetermined
number, $K_j$, of Markov transitions at stage $j$, we might instead
perform as many transitions as are necessary to obtain a good sample.
Define a `tour' to be a sequence of transitions that moves from a high
value of some key quantity (eg, $U(x)$ for the canonical distributions
of equation~(\ref{eq-canonical})) to a low value of this quantity, or
vice versa.  Good sampling might be ensured by performing some
predetermined number of tours, with the number of these tours that
occur before and after the link state being chosen at random.
Suitable `high' and `low' values would probably need to be found using
preliminary runs.

More speculatively, it seems as if there should be some method that
has the advantages of LIS over AIS, but that like AIS uses many
intermediate distributions, performing only a single Markov transition
for each.  Intuitively, it seems that such a `smooth' method that does
not abruptly change $\eta$ should be more efficient.  One can use LIS
with all $K_j$ set to one, but this will produce good results only if
$n$ is large, which we saw in the analysis of Section~\ref{sec-asym}
does not lead to an advantage over AIS.  Perhaps some way could be
found of using states associated with all values of $\eta$ when
estimating each of the ratios $Z_{\eta_{j+1}}/Z_{\eta_j}$, while still
producing an estimate that is exactly unbiased even when the Markov transitions
do not reach equilibrium.

\section*{Acknowledgements}\vspace{-10pt}

This research was supported by the Natural Sciences and Engineering
Research Council of Canada.  I hold a Canada Research Chair in
Statistics and Machine Learning.

\section*{References}\vspace{-10pt}

\leftmargini 0.2in
\labelsep 0in

\begin{description}
\itemsep 2pt

\item
  Bennett, C.~H.\ (1976) ``Efficient estimation of free energy differences
  from Monte Carlo data'', {\em Journal of Computational Physics}, vol.~22,
  pp.~245-268.

\item 
  Crooks, G.~E.\ (2000) ``Path-ensemble averages in systems driven far
  from equilibrium'', \textit{Physical Review E}, vol.~61, pp.~2361-2366.


\item
  Gelman, A.\ and Meng, X.-L.\ (1998) ``Simulating normalizing constants:
  From importance sampling to bridge sampling to path sampling'', 
  \textit{Statistical Science}, vol.~13, pp.~163-185.

\item
  Hendrix, D.~A.\ and Jarzynski, C.\ (2001) ``A ``fast growth'' method of
  computing free energy differences'', \textit{Journal of Chemical Physics},
  vol.~114, pp.~5974-5981.

\item
  Jarzynski, C.\ (1997) ``Nonequilibrium equality for free energy differences'',
  \textit{Physical Review Letters}, vol.~78, pp.~2690-2693.

\item
  Jarzynski, C.\ (2001) ``A ``fast growth'' method of computing free energy
  differences'', \textit{Journal of Chemical Physics}, vol.~114, pp.~5974-5981.


\item
  Lu, N., Singh, J.~K., and Kofke, D.~A.\ (2003) ``Appropriate methods
  to combine forward and reverse free-energy perturbation averages'',
  \textit{Journal of Chemical Physics}, vol.~118, pp.~2977-2984.

\item
  MacEachern, S.~N., Clyde, M., and Liu, J.~S. (1999) ``Sequential
  importance sampling for nonparametric Bayes models:\ The next generation'',
  \textit{Canadian Journal of Statistics}, vol.~27, pp.~251-267.

\item
  Meng, X.-L.\ and Wong, H.~W.\ (1996) ``Simulating ratios of normalizing
  constants via a simple identity: A theoretical exploration'', 
  \textit{Statistica Sinica}, vol.~6, pp.~831-860.

\item
  Metropolis, N., Rosenbluth, A.~W., Rosenbluth, M.~N., Teller, A.~H., 
  and Teller, E.\ (1953) ``Equation of state calculations by fast computing 
  machines'', {\em Journal of Chemical Physics}, vol.~21, pp.~1087-1092.

\item
  Neal, R.~M.\ (1993) {\em Probabilistic Inference Using Markov Chain
  Monte Carlo Methods}, Technical Report CRG-TR-93-1, Dept.\
  of Computer Science, University of Toronto, 140 pages.  
  Obtainable from \texttt{http://www.cs.utoronto.ca/$\sim$radford/}.

\item
  Neal, R.~M.\ (1996) ``Sampling from multimodal distributions using tempered 
  transitions'', \textit{Statistics and Computing}, vol.~6, pp.~353-366.

\item
  Neal, R.~M.\ (1999) ``Regression and classification using Gaussian process
      priors'' (with discussion), in J.~M.~Bernardo, {\em et al} 
      (editors) {\em Bayesian Statistics 6}, Oxford University Press, 
      pp.~475-501.

\item
  Neal, R.~M.\ (2001) ``Annealed importance sampling'', 
  \textit{Statistics and Computing}, vol.~11, pp.~125-139.


\item
  Neal, R.~M.\ (2004) ``Taking bigger Metropolis steps by dragging fast 
  variables'', Technical Report No.~0411, Dept. of Statistics, University of 
  Toronto, 9 pages.

\item
  Schervish, M.~J.\ (1995) \textit{Theory of Statistics}, Springer.

\item 
  Shirts, M.~R., Bair, E., Hooker, G., and Pande, V.~S.` (2003)
  ``Equilibrium free energies from nonequilibrium measurements using
  maximum-likelihood methods'', \textit{Physical Review Letters},
  vol.~91, p.~140601.

\end{description}

\end{document}